# Localization and delocalization of random interfaces[*]

## Yvan Velenik

*CNRS and Université de Rouen*

*Avenue de l'Université, BP.12*
*76801 Saint Étienne du Rouvray (France)*
*e-mail:* Yvan.Velenik@univ-rouen.fr
*url:* www.univ-rouen.fr/LMRS/Persopage/Velenik/

**Abstract:** The probabilistic study of effective interface models has been quite active in recent years, with a particular emphasis on the effect of various external potentials (wall, pinning potential, ...) leading to localization/delocalization transitions. I review some of the results that have been obtained. In particular, I discuss pinning by a local potential, entropic repulsion and the (pre)wetting transition, both for models with continuous and discrete heights.

**AMS 2000 subject classifications:** Primary 60K35, 82B41; secondary 82B24, 82B27.

**Keywords and phrases:** random interfaces, random surfaces, roughening, entropic repulsion, pinning, wetting, prewetting.

## Contents

---

[*]These notes were prepared for a 6 hours mini-course given at the workshop *Topics in Random Interfaces and Directed Polymers*, Leipzig, September 12-17 2005.

## 1. Introduction

In statistical physics, macroscopic systems in thermodynamic equilibrium are described by a special type of random fields, the Gibbs measures, that are defined through a family of potentials, encoding the physics of the system at the microscopic scale [58]. It may turn out that a given family of potentials gives rise to several distinct Gibbs states; in that case, one says that there is phase coexistence. This corresponds to the well-known phenomenon of first-order phase transitions, when several phases are thermodynamically stable at the same value of the relevant thermodynamical parameters (temperature, pressure, etc.): for example (under standard conditions) water and ice are both stable at 0°C, while liquid water and water vapor are both stable at 100°C. A particularly interesting feature of systems in the regime of phase coexistence is that it is possible to have different regions of space occupied by different phases, and separated from each other by sharp boundaries, the interfaces. The latter are the central objects of these notes.

The most natural way to analyze such interfaces is to consider a lattice gas in the regime of phase coexistence; for example, a $d$-dimensional, $d \geq 2$, Ising model below the critical temperature. In such models, one can then partition space into cells of suitable size, and attribute a label to each cell, depending on the behaviour inside the cell. For example, in the Ising model, one might get three types of labels: 1 when the system in the cell "looks like" the positively magnetized state, $-1$ when it "looks like" the negatively magnetized state, and 0 otherwise. If the size of the cells is chosen carefully, the connected components of cells with 0-label may provide a possible definition of interfaces. We see that in these models interfaces are emergent quantities appearing because the spatial coexistence of several Gibbs states is enforced by suitable constraints (boundary conditions, fixed magnetization, etc.). The study of interfaces is thus impaired by the fact that one has first to locate them (and this can only be done with ambiguity), and then the resulting objects will have a complicated structure (finite but possibly large microscopic thickness with non-trival internal structure).

Surprisingly enough, as long as one is interested in the *macroscopic* description of the system, the analysis can be pushed quite far. In particular, it is possible to provide a precise description of the macroscopic geometry of these interfaces, and in particular to prove that the latter is usually deterministic and given by the solution of suitable variational problems involving the relevant surface free energies. We refer to the review paper [12] for much more information on this topic.

The situation is substantially less satisfactory when one is interested in meso- or microscopic properties of the interfaces, such as the statistical description of their fluctuations. Even though it is still possible to obtain precise informations in two dimensions (at least in several interesting cases), the situation is almost completely out of reach in higher dimensions (except in some very special situations at very low temperatures). In order to make progress on these issues it is thus useful to consider a class of more tractable models, in which the interfaces are not emergent quantities anymore, but are the basic microscopic objects one works with. A natural class, and by far the most studied, of such *effective interface models* is the one discussed in these notes.

In this class of models, the interface is modeled as the graph of a random function from $\mathbb{Z}^d$ to $\mathbb{Z}$ or $\mathbb{R}$ (discrete and continuous effective interface models), and its distribution is given by a Gibbs measure with an interaction depending only on the gradients of the field. These models are expected to provide good approximations of "real" interfaces, at least in suitable regimes. Namely, computations with discrete effective interface models can be shown to provide good agreement with interfaces of lattice gases at low temperature as long as both models are below their roughening transition (see Subsection 1.3). "Real" inter-faces above their roughening transition should be well described by continuous effective interface models, but we lack rigorous proof of that, except for some very special cases (see Subsection 1.3). The only situation where this can be sometimes proved is $d = 1$. In that case, it is for example possible to prove that interfaces in the 2D Ising model have the same Brownian asymptotics as their effective counterparts [66].

This type of models is also used by physicists to make predictions (usually based on exact computations for one-dimensional interfaces, or non-rigorous renormalization group methods for higher dimensional interfaces) on the behaviour of interfaces in real, or more realistic, systems. This works quite well in general, even though there are problems when investigating some phenomena (in particular, the wetting transition, see Section 6), and it turned out to be necessary to introduce more complicated models.

There are also good reasons for a probabilist to be interested in such models. First, they correspond to Gibbs random fields with truly unbounded spins, and the latter property gives rise to many phenomena not encountered when studying finite or compact spin-spaces. Second, these models of random surfaces can be seen as natural generalizations of random walks to higher-dimensional "time". Indeed, the Gibbs measures of one-dimensional effective interface models with nearest-neighbor interactions are nothing else but the path measures of random walks (usually with bridge-like conditionning, since both endpoints are fixed): After all, the path of a discrete-time random walk on $\mathbb{Z}$ (or $\mathbb{R}$) is the graph of a random function from $\mathbb{Z}$ to $\mathbb{Z}$ (or $\mathbb{R}$), and the weight of a trajectory obviously only depends on the increments (that is, the gradients of the path).

Many interesting results have been obtained for these effective interface models during the last decade. It is not possible to discuss all the investigated topics in these notes, so here the main focus is on the localization/delocalization transitions for these interfaces. Namely, after defining and stating the basic properties of such random surfaces, I'll discuss the effect of various external potentials of physical relevance on their behaviour. Often, phase transitions between localized and delocalized states occur, and understanding and quantifying them is the main theme of these notes.

Fortunately, there are many good references for people interested in other topics. The derivation of the variational problems describing the shape of interfaces at macroscopic scale (Wulff shape, etc.) is discussed in the review paper [12] (mostly for lattice gases) and the Saint-Flour lecture notes of Funaki [55]. The dynamical aspects of these models are of great interest and are also discussed in [55]. Moreover, there is a very close connection between models of polymer chains and one-dimensional effective interface models; an excellent reference on polymers for probabilists is the book by Giacomin [62].

Finally, concerning the topics discussed in these notes, additional informations on continuous effective interface models can be found in the lecture notes of Giacomin [60], in the last part of the Saint-Flour lecture notes of Bolthausen [14], and in his enlightening review on entropic repulsion [13]. There are not nearly as many references for discrete effective interface models, but the older and excellent reviews by Fisher [50] and Bricmont *et al* [26] still provide a lot of informations.

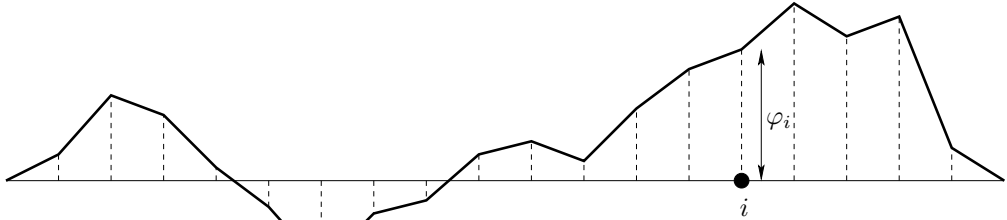

FIG 1. *A one-dimensional interface $\varphi$.*

### 1.1. The free model

The interface is described by a function $\varphi : \mathbb{Z}^d \to \mathbb{R}$ (*continuous* effective interface models) or $\varphi : \mathbb{Z}^d \to \mathbb{Z}$ (*discrete* effective interface models); the quantity $\varphi_i$ is interpreted as the height of the interface above (or below) site $i$; see Fig 1. To each realization of the random field and any $\Lambda \Subset \mathbb{Z}^d$, we associate an energy given by the Hamiltonian

$$H_\Lambda(\varphi) = \tfrac{1}{2} \sum_{i,j \in \Lambda} p(j-i) V(\varphi_j - \varphi_i) + \sum_{i \in \Lambda, j \notin \Lambda} p(j-i) V(\varphi_j - \varphi_i), \quad (1)$$

Here:

- $V : \mathbb{R} \to \mathbb{R}$ is an even, convex function, satisfying $V(0) = 0$. Often, it will be assumed that it is actually *uniformly strictly convex*, which means that $V$ is $C^2$ and $c_- < V''(x) < c_+$ for some $0 < c_- < c_+ < \infty$.
- $p(\,\cdot\,)$ is the transition kernel of an (usually aperiodic, irreducible) symmetric discrete-time random walk $X$ on $\mathbb{Z}^d$; the law of $X$ started from site $i$ at time 0 is denoted by $\mathbb{P}_i$ and its expectation $\mathbb{E}_i$. We'll say that the model has *finite-range* interactions if $X$ has bounded jumps; we'll say that it has *nearest-neighbor* interactions if $p(i) = 0$ when $\|i\|_1 \neq 1$. When the range of the interaction is not finite, we'll always assume that there exists $\delta > 0$ such that

$$\mathbb{E}_0 \left( \|X_1\|_2^{d+\delta} \right) < \infty. \quad (2)$$

Sometimes, we'll have to make the stronger assumption that there exists $\alpha > 0$ such that

$$\mathbb{E}_0 \left( e^{\alpha \|X_1\|_2} \right) < \infty. \quad (3)$$

We'll use the notation $\mathcal{Q}$ for the $d \times d$ matrix with elements $\mathcal{Q}_{rs} = \mathbb{E}_0 \left( X_1(r) X_1(s) \right)$, where $x(r)$ is the $r$-th component of $x \in \mathbb{R}^d$.

The statistical properties of the free interface are described by the following Gibbs measure on $\mathbb{R}^{\mathbb{Z}^d}$, resp. $\mathbb{Z}^{\mathbb{Z}^d}$. Let $\Lambda \Subset \mathbb{Z}^d$ (this notation meaning that $\Lambda$ is a finite subset of $\mathbb{Z}^d$), let $\psi$ be a fixed realization of the field, and let $\beta \geq 0$. In the case of continuous heights, the corresponding finite-volume Gibbs measure is defined as the following probability measure on $\mathbb{R}^{\mathbb{Z}^d}$ (with the product Borel $\sigma$-algebra)

$$\mathrm{P}_\Lambda^{\psi,\beta}(\mathrm{d}\varphi) \stackrel{\text{def}}{=} \left(Z_\Lambda^{\psi,\beta}\right)^{-1} \exp\left(-\beta H_\Lambda(\varphi)\right) \prod_{i \in \Lambda} \mathrm{d}\varphi_i \prod_{j \notin \Lambda} \delta_{\psi_j}(\mathrm{d}\varphi_j), \qquad (4)$$

where $\delta_x$ denotes Dirac's mass at $x$, and the *partition function* $Z_\Lambda^{\psi,\beta}$ is the normalization constant. The corresponding probability measure in the case of discrete heights is then given by

$$\mathrm{P}_\Lambda^{\psi,\beta}(\varphi) \stackrel{\text{def}}{=} \left(Z_\Lambda^{\psi,\beta}\right)^{-1} \exp\left(-\beta H_\Lambda(\varphi)\right) \prod_{j \notin \Lambda} \delta_{\psi_j,\varphi_j},$$

where $\delta_{i,j}$ is Kronecker's $\delta$-function. Expectation with respect to $\mathrm{P}_\Lambda^{\psi,\beta}$ is denoted by $\mathrm{E}_\Lambda^{\psi,\beta}$.

Apart from some remarks, these lectures are restricted to a few important special cases: In the continuous case, our main example will be $V(x) = \frac{1}{2}x^2$, to which we refer as the case of quadratic, or Gaussian, interactions; notice that in this case, the measure $\mathrm{P}_\Lambda^{\psi,\beta}$ is actually Gaussian, which makes this model particularly tractable. In the discrete case, we will focus on $V(x) = |x|$, the *Solid-On-Solid* (SOS) model, or $V(x) = \frac{1}{2}x^2$, the *discrete Gaussian* (DG) model.

**Remark 1.** *In the (continuous) Gaussian case, the inverse temperature $\beta$ can be set to 1 by a simple change of variables. Therefore in the sequel, when stating results about this model, it will always be assumed that this has been done, and $\beta$ will be removed from the notation.*

**Remark 2.** *In these notes, I mostly consider 0 boundary conditions, $\psi \equiv 0$. In that case, I will therefore omit any mention of the boundary conditions in the notations.*

### 1.2. Continuous effective interface models: basic properties

I start by discussing basic properties of continuous effective interface models, as these are in general better understood than their discrete counterparts, and should be closely related to the latter in important cases. I restrict the discussion to the case of 0-boundary conditions, $\psi \equiv 0$. I also only consider the Gaussian case, except for a few remarks.

*1.2.1. The random walk representation*

As pointed out before, in the case of quadratic $V$, the measure $P_\Lambda$ is Gaussian. As such, it is characterized by its covariance matrix. It turns out that there is a very nice, and extremely useful representation of the latter in terms of the Green function of the random walk $X$ with transition kernel $p(\,\cdot\,)$.

Elementary algebraic manipulations yield

$$H_\Lambda(\varphi) = \tfrac{1}{4} \sum_{i,j\in\Lambda} p(j-i)(\varphi_j - \varphi_i)^2 + \tfrac{1}{2} \sum_{i\in\Lambda, j\notin\Lambda} p(j-i)(\varphi_j - \varphi_i)^2$$
$$= \tfrac{1}{2} \sum_{i,j\in\Lambda} \varphi_i(\delta_{ij} - p(j-i))\varphi_j\,.$$

where $(\mathbf{1} - \mathbf{P})_{i,j} \stackrel{\text{def}}{=} \delta_{i,j} - p(j-i)$; Consequently, introducing the matrices $\mathbf{1}_\Lambda \stackrel{\text{def}}{=} (1)_{i,j\in\Lambda}$ and $\mathbf{P}_\Lambda \stackrel{\text{def}}{=} (p(j-i))_{i,j\in\Lambda}$, the covariances are seen to be given by

$$\operatorname{cov}_{P_\Lambda}(\varphi_i, \varphi_j) = \left[(\mathbf{1}_\Lambda - \mathbf{P}_\Lambda)^{-1}\right]_{ij} = \left[\sum_{n\geq 0} \mathbf{P}_\Lambda^n\right]_{ij}$$
$$= \sum_{n\geq 0} \mathbb{P}_i[X_n = j, \tau_\Lambda > n] \equiv G_\Lambda(i,j)\,, \quad (5)$$

where $X$ is the random walk with transition kernel $p$ introduced above, and $\tau_\Lambda \stackrel{\text{def}}{=} \min\{n : X_n \notin \Lambda\}$. $G_\Lambda$ is thus the Green function of this random walk, killed as it exits $\Lambda$.

**Remark 3.** *Similar representations are available for the mean $E_\Lambda^\psi(\varphi_i)$, and for the partition function $Z_\Lambda^\psi$.*

**Remark 4.** *This derivation of the random walk representation of course relies entirely on the Gaussian nature of the measure. Nevertheless, it turns out that it is possible to derive a generalization of this representation valid for the whole class of uniformly strictly convex interactions. Such a generalization has been proposed in [41] and is a probabilistic reformulation of an earlier result, in the PDE context, by Helffer and Sjöstrand [67]. It works as follows: One constructs a stochastic process $(\Phi(t), X(t))$ where*

- *$\Phi(\,\cdot\,)$ is a diffusion on $\mathbb{R}^{\mathbb{Z}^d}$ with invariant measure $P_\Lambda$;*
- *given a trajectory $\varphi(\,\cdot\,)$ of the process $\Phi$, $X(t)$ is an, in general inhomogeneous, transient, continuous-time random walk on $\mathbb{Z}^d$ with life-time $\tau_\Lambda \stackrel{\text{def}}{=} \inf\{t \geq 0 \,|\, X(t) \notin \Lambda\}$, and time-dependent symmetric jump-rates*

$$a(i,j;t) \stackrel{\text{def}}{=} p(j-i)V''(\varphi_j(t) - \varphi_i(t))\,.$$

*Denoting by $\mathcal{E}_{i,\varphi}^\Lambda$ the law of $(\Phi(t), X(t))$ starting from the point $(i, \varphi) \in \Lambda \times \mathbb{R}^{\mathbb{Z}^d}$, we have the following generalization of* (5),

$$\operatorname{cov}_\Lambda(\varphi_i, \varphi_j) = E_\Lambda\left(\mathcal{E}_{i,\cdot}^\Lambda \int_0^{\tau_\Lambda} \mathbf{1}_{\{X(s)=j\}}\mathrm{d}s\right). \quad (6)$$

*Thanks to the ellipticity of the random walk $X(t)$ under the assumption of uniform strict convexity, it is possible to obtain some Aronson type bounds, see [63, 37], showing that this random walk in a (dynamical) random environment has the same qualitative behavior as the random walk in the Gaussian case. This explains why most of the results that have been obtained for the Gaussian model also hold in the non-Gaussian case. However, since quantitative estimates are still out of reach for this process, the quantitative Gaussian results are replaced by qualitative ones.*

### 1.2.2. Transverse and longitudinal correlation lengths

Physicists often characterize the statistical properties of interfaces through two quantities: the *transverse correlation length*, which measures the fluctuation in the direction orthogonal to the surface, and the *longitudinal correlation length*, which measures the correlations along the surface. From the mathematical point of view, these two quantities are directly related to, respectively, the variance of the field, and the rate of exponential decay of its covariance, also called in the physics literature the *mass*. The mass associated to an infinite-volume Gibbs measure Q is defined, for any $x \in \mathbb{S}^{d-1}$ by

$$m_Q(x) \stackrel{\text{def}}{=} - \lim_{k \to \infty} \frac{1}{k} \log \text{cov}_Q(\varphi_0, \varphi_{[kx]}), \tag{7}$$

In the last expression, we have used the notation $[x] \stackrel{\text{def}}{=} ([x(1)], \ldots, [x(d)]) \in \mathbb{Z}^d$ when $x = (x(1), \ldots, x(d)) \in \mathbb{R}^d$. Of course, it is not obvious that this limit exists in general. However in most cases of interest, this can be proved using suitable subadditivity arguments. When this is not the case, the bounds on $m$ given in the sequel have to be understood as follows: Lower bounds on $m$ must be interpreted as upper bounds on the quantity obtained by replacing lim by lim sup, while upper bounds on $m$ must be interpreted as upper bounds on the quantity obtained by replacing lim by lim inf.

Let us see how these quantities behave in the case of the free interface. Let $\Lambda_N = \{-N, \ldots, N\}^d$. It follows from the random walk representation (5) and standard estimates for the Green function of the corresponding random walk [88, 76] that there exist constants $0 < g_d < \infty$, depending on the transition kernel $p(\cdot)$, such that, as $N \to \infty$,

$$\text{var}_{P_{\Lambda_N}}(\varphi_0) = \begin{cases} (g_1 + o(1))\, N & (d = 1), \\ (g_2 + o(1)) \log N & (d = 2), \\ g_d + o(1) & (d \geq 3). \end{cases} \tag{8}$$

We see that the variance diverges when $d = 1$ or $2$ (the corresponding interface is said to be *delocalized*), while it remains finite for $d \geq 3$ (the interface is localized). It also follows from these results that the limiting field P exists if and only if $d \geq 3$.

When $d \geq 3$, the random walk representation also provides informations on the covariances of the infinite-volume field; namely [88], there exists $R_d > 0$, also depending on the transition kernel $p(\,\cdot\,)$, such that, as $|i| \to \infty$,

$$\mathrm{cov}_{\mathrm{P}}(\varphi_0, \varphi_i) = (R_d + o(1)) \, \|i\|_{\mathcal{Q}}^{2-d},$$

where we used the norm $\|i\|_{\mathcal{Q}} = (i, \mathcal{Q}i)^{1/2}$. (This result requires that the transition kernel has slightly more than moments of order $d$ [76], a condition satisfied when (3) holds.) We see that the corresponding limiting field has very strong correlations (not summable). In particular, the mass satisfies $m_{\mathrm{P}}(x) \equiv 0$. This is the reason why these models are usually called *massless*.

To summarize, the free interface is delocalized when $d = 1$ or $2$, and localized but strongly correlated when $d \geq 3$.

**Remark 5.** *All the results in this subsection also hold true in the general case of uniformly strictly convex interactions, at least qualitatively (i.e. one has upper and lower bounds for these quantities that are of the same order as those in the Gaussian settings, but which do not coincide). Actually, in the Gaussian case, much more precise results than those given here can be obtained.*

**Remark 6.** *Apart from the random walk representation, I am aware of only one general tool to prove localization: the Brascamp-Lieb inequality, see [25]. Unfortunately, the class to which this approach applies, if already quite large, is still much too limited. Namely, it is required that $V$ satisfies one of the following conditions:*

1. *$V(x) = ax^2 + f(x)$, with $f$ convex and $a > 0$;*
2. *$0 < c_- \leq V''(x) \leq c_+ < \infty$ for all $|x| \geq M$, for some $M < \infty$, and $|V(x) - Cx^2| \leq D < \infty$, for all $x$ and some $C > 0$.*

*Concerning the proof of delocalization, on the other hand, several methods exist, all with their own limitations. One of them is the subject of the next subsection.*

**Open Problem 1.** *Prove localization in high dimensions for a more general class of interactions $V$ than those treated in [25]. For example, the case $V(x) = x^4$ is still open.*

*1.2.3. Delocalization as a consequence of continuous symmetry*

There is another point of view, which is interesting in order to understand the delocalization of the interface in dimensions 1 and 2, and which shows that this should be a very common phenomenon, valid for very general interactions $V$. The main observation is that the Hamiltonian enjoys a continuous symmetry: $H_\Lambda(\varphi) = H_\Lambda(\varphi + c)$, for any $\Lambda \Subset \mathbb{Z}^d$ and $c \in \mathbb{R}$, since it is actually only a function of the gradient field $\nabla\varphi$. By standard Mermin-Wagner–type arguments [45, 83, 25, 69], it then follows that this continuous symmetry should also be present at the level of the infinite-volume Gibbs measures, when $d = 1$ or $2$ and the interaction does not decay too slowly (a condition automatically satisfied under (2)). Of course, this is impossible, since it would imply, for example, that the law of $\varphi_0$ under this measure is uniform on $\mathbb{R}$; this means that there cannot be any infinite-volume Gibbs measure when $d = 1$ or $2$. Actually, it is even possible to derive, using such arguments, qualitatively correct lower bounds on the size of the fluctuations of the interface. I show here how such a claim is proved when $V$ is twice continuously differentiable and such that $V'' \leq c$. To see how to treat more general $V$ and general boundary conditions, I refer to [69].

**Remark 7.** *It should be emphasized, however, that even though a large class of interactions $V$ can be treated in this way, much too strong assumptions on $V$ are still required, and as such the situation is still far from satisfactory. The alternative approach, pioneered in [25] and based on Brascamp-Lieb inequalities, also imposes unsatisfactory (though different) assumptions on $V$. More precisely, current methods of proof require that $V$ satisfies one of the following conditions:*

1. *$\|V - \bar{V}\|_\infty < \epsilon$ for some small enough $\epsilon$ and some twice continuously differentiable function $\bar{V}$ such that $\bar{V}''(x) \leq c < \infty$, for all $x$;*
2. *$\lim_{|x|\to\infty} (|x| + |V'(x)|) \exp(-V(x)) = 0$, and $|V'(x)| \leq c < \infty$, for all $x$;*
3. *$\lim_{|x|\to\infty} (|x| + |V'(x)|) \exp(-V(x)) = 0$, $V$ is convex and*

$$\int (V'(x))^2 \exp(-V(x)) < \infty.$$

**Open Problem 2.** *Prove delocalization in low dimensions for a larger class of interactions than those treated in [25] and [69]. For example, the interesting case of the Hammock potential, $V(x) = 0$ for $|x| < 1$ and $V(x) = \infty$ for $|x| > 1$, is still open.*

Let us fix some configuration $\bar{\varphi}$, such that $\bar{\varphi}_i = 0$, for all $i \notin \Lambda_N$, and $\bar{\varphi}_0 = R$, for some $R > 0$. We'll choose $R$ and optimize over $\bar{\varphi}$ later on. Let us introduce the transformation $\mathsf{T}_{\bar{\varphi}}(\varphi) \stackrel{\text{def}}{=} \varphi + \bar{\varphi}$, and the tilted measure $\mathrm{P}_{\Lambda_N;\bar{\varphi}} \stackrel{\text{def}}{=} \mathrm{P}_{\Lambda_N} \circ \mathsf{T}_{\bar{\varphi}}$. Observe that

$$\mathrm{P}_{\Lambda_N}(\varphi_0 \geq R) = \mathrm{P}_{\Lambda_N;\bar{\varphi}}(\varphi_0 \geq 0). \tag{9}$$

Recall the standard entropy inequality (actually a simple consequence of Jensen inequality, see, *e.g.*, [60, Appendix B.3])

$$\mu(A) \geq \nu(A) \exp\left(-(\mathbf{H}(\nu|\mu) + e^{-1})/\nu(A)\right), \tag{10}$$

where $\mu$ and $\nu$ are probability measures such that $\nu \ll \mu$, $A$ is some event with $\nu(A) > 0$, and $\mathbf{H}(\nu|\mu) \stackrel{\text{def}}{=} \mathrm{E}_\nu(\log \frac{\mathrm{d}\nu}{\mathrm{d}\mu})$ is the relative entropy of $\nu$ w.r.t. $\mu$. Applying this inequality, we get that

$$\mathrm{P}_{\Lambda_N;\bar{\varphi}}(\varphi_0 \geq 0) \geq \mathrm{P}_{\Lambda_N}(\varphi_0 \geq 0) \exp\left(-\left(\mathbf{H}(\mathrm{P}_{\Lambda_N;\bar{\varphi}}|\mathrm{P}_{\Lambda_N}) + e^{-1}\right)/\mathrm{P}_{\Lambda_N}(\varphi_0 \geq 0)\right)$$
$$= \tfrac{1}{2} \exp\left(-2(\mathbf{H}(\mathrm{P}_{\Lambda_N;\bar{\varphi}}|\mathrm{P}_{\Lambda_N}) + e^{-1})\right).$$

Therefore, it follows from (9) that

$$\mathrm{P}_{\Lambda_N}(\varphi_0 \geq R) \geq \tfrac{1}{2} \exp\left(-2(\mathbf{H}(\mathrm{P}_{\Lambda_N;\bar{\varphi}}|\mathrm{P}_{\Lambda_N}) + e^{-1})\right),$$

and it only remains to control the relative entropy. A simple computation yields (observe that the partition functions are equal, since the Jacobian of the transformation is 1)

$$\begin{aligned}\mathbf{H}(\mathrm{P}_{\Lambda_N;\bar{\varphi}}|\mathrm{P}_{\Lambda_N}) &= \mathrm{E}_{\Lambda_N}\left(H(\varphi - \bar{\varphi}) - H(\varphi)\right) \\ &= \tfrac{1}{2}\sum_{i,j\in\Lambda} p(j-i)\, \mathrm{E}_{\Lambda_N}\left(V(\varphi_j - \varphi_i - \bar{\varphi}_j + \bar{\varphi}_i) - V(\varphi_j - \varphi_i)\right) + \\ &\quad + \sum_{i\in\Lambda, j\notin\Lambda} p(j-i)\, \mathrm{E}_{\Lambda_N}\left(V(\varphi_j - \varphi_i - \bar{\varphi}_j + \bar{\varphi}_i) - V(\varphi_j - \varphi_i)\right) \\ &\leq \frac{c}{4}\sum_{i,j\in\Lambda} p(j-i)\,(\bar{\varphi}_j - \bar{\varphi}_i)^2 + \frac{c}{2}\sum_{i\in\Lambda, j\notin\Lambda} p(j-i)\,(\bar{\varphi}_j - \bar{\varphi}_i)^2\,.\end{aligned}$$

The last inequality follows from a Taylor expansion, keeping in mind that $\mathrm{E}_{\Lambda_N}(\varphi_j - \varphi_i) = 0$ by symmetry and that $V'' \leq c$.

Therefore, choosing $\bar{\varphi}_i \stackrel{\mathrm{def}}{=} R\,\mathbb{P}_i\left[\mathrm{T}_{\{0\}} < \tau_{\Lambda_N}\right]$, where $\mathrm{T}_{\{0\}} \stackrel{\mathrm{def}}{=} \min\{n\,:\, X_n = 0\}$, and using the following well-known estimate [76] for two-dimensional random walks with finite variance[1]

$$\mathbb{P}_i\left[\mathrm{T}_{\{0\}} > \tau_{\Lambda_N}\right] \asymp \frac{\log(|i|+1)}{\log(N+1)}\,,$$

we see that $\mathbf{H}(\mathrm{P}_{\Lambda_N;\bar{\varphi}}|\mathrm{P}_{\Lambda_N}) = O(R^2/\log N)$, and therefore

$$\mathrm{P}_{\Lambda_N}(\varphi_0 > T\sqrt{\log N}) \geq \exp(-cT^2)\,.$$

### 1.2.4. Tail of one-site marginals

When the interface is localized, the finiteness of the variance provides only rather weak information on the fluctuations. Another quantity that is of interest is the tail of the one-site marginal, as this shows how strongly the interface is localized. As we have just seen, the free interface is only localized when $d \geq 3$. Evidently, when the infinite-volume measure is Gaussian, $\varphi_0$ is just a Gaussian random variable with variance $\sum_{n\geq 0}\mathbb{P}_0(X_n = 0) < \infty$, and therefore its has a Gaussian tail.

**Remark 8.** *Of course, the tails have also Gaussian decay in the case of uniformly strictly convex interaction $V$, as follows, e.g., from Brascamp-Lieb inequality [24].*

### 1.2.5. Thermodynamical criteria of localization

The above quantities are those one would very much like to compute in every situations. Unfortunately, this often turns out to be too difficult, and we have to content ourselves with much more limited informations, often at the level of the free energy. These thermodynamical criteria work by comparing the free energy of the system under consideration with that of the free interface. They can usually be reinterpreted (by the usual "differentiate, then integrate back" technique) as estimates on the expected value of a suitable macroscopic quantity. Let us consider a simple example of the latter.

---

[1]The notation $a \asymp b$ means that there exist two constants $0 < c_1 \leq c_2 < \infty$, depending only on the dimension, such that $c_1 b \leq a \leq c_2 b$.

Let $\Delta > 0$, and set $\rho_N \stackrel{\text{def}}{=} |\Lambda_N|^{-1} \sum_{i \in \Lambda_N} \mathbf{1}_{\{|\varphi_i| \leq \Delta\}}$. One could then say that the interface is localized when $\liminf_N \rho_N > 0$ and delocalized otherwise. By itself, such an estimate does not guarantee that the interface is really localized or delocalized, in a pathwise sense, but it is a strong indication that it should be the case. Note that for the free interface, this yields again localization if and only if $d \geq 3$.

A similar thermodynamical criterion is used when discussing the wetting transition in Section 6.

### 1.2.6. Some additional properties of the Gaussian model

**Extrema of the field.** In order to understand properly the entropic repulsion phenomenon discussed in Section 3, it is important to know the behavior of the large fluctuations of the interface. These have a very different behavior in dimension 1, 2, and in higher dimensions.

One is not really interested in the behavior when $d = 1$, as the latter is dictated by the behavior of the corresponding quantity for the Brownian bridge. In higher dimension the results are more interesting, since the extrema of the field turn out to be much larger than the typical values. It is proved in [16] that the maximum of the 2-dimensional finite-range Gaussian field in the box $\Lambda_N$ satisfies, for any $\delta > 0$,

$$\lim_{N \to \infty} \mathrm{P}_{\Lambda_N} \left( \left| \sup_{i \in \Lambda_N} \varphi_i - 2\sqrt{g_2} \log N \right| \geq \delta \log N \right) = 0, \tag{11}$$

where $g_d$ was introduced in (8). Similarly, it is proved in [17, 18] that the maximum of the $d$-dimensional, $d \geq 3$, finite-range Gaussian field in the box $\Lambda_N$ satisfies, for any $\delta > 0$,

$$\lim_{N \to \infty} \mathrm{P}_{\Lambda_N} \left( \left| \sup_{i \in \Lambda_N} \varphi_i - \sqrt{2dg_d}\sqrt{\log N} \right| \geq \delta \sqrt{\log N} \right) = 0. \tag{12}$$

It is interesting to note that one bound is actually obvious, since, *e.g.*,

$$\mathrm{P}_{\Lambda_N} \left( \sup_{i \in \Lambda_N} \varphi_i > (2\sqrt{g_2} + \delta) \log N \right) \leq |\Lambda_N| \sup_{i \in \Lambda_N} \mathrm{P}_{\Lambda_N} \left( \varphi_i > (2\sqrt{g_2} + \delta) \log N \right),$$

which vanishes as $N \to \infty$ (use the fact that $\mathrm{var}_{\mathrm{P}_{\Lambda_N}}(\varphi_i) \leq \mathrm{var}_{\mathrm{P}_{\Lambda_N}}(\varphi_0)$).

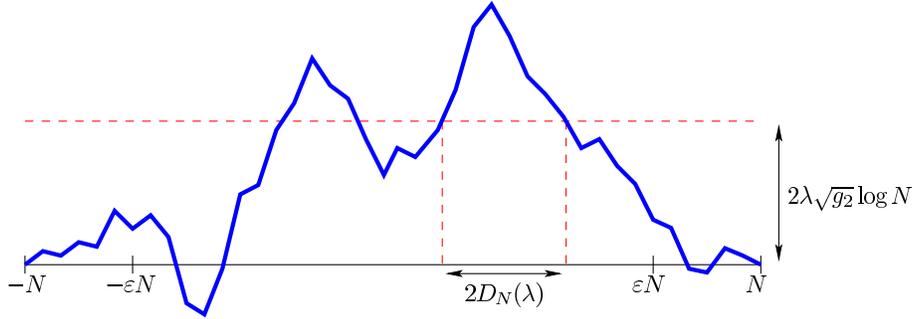

FIG 2. *$D_N(\lambda)$ measures the width of the widest spike reaching above $2\lambda\sqrt{g_2}\log N$ in the subbox $\Lambda_{\varepsilon N}$ (of course, the real situation is 2-dimensional).*

**The spikes.** We have just seen that the $d$-dimensional interface ($d \geq 2$) can make huge fluctuations, called spikes. It is interesting to investigate the geometry and distribution of these spikes.

The most interesting picture emerges when $d = 2$. A detailed study of the spikes of the 2-dimensional nearest-neighbor Gaussian model can be found in [34]. The main results can be stated as follows.

- The spikes are rather fat. Let $0 < \lambda < 1$, $0 < \varepsilon < 1$, and (see Fig. 2)

$$D_N(\lambda) \stackrel{\text{def}}{=} \sup\left\{ a \in \mathbb{N} \,:\, \exists i \in \Lambda_{\varepsilon N},\, \min_{\|j-i\|_\infty \leq a} \varphi_j > 2\lambda\sqrt{g_2}\log N \right\}.$$

Then
$$\lim_{N\to\infty} \frac{\log D_N(\lambda)}{\log N} = \frac{1}{2} - \frac{\lambda}{2} \qquad \text{in probability.}$$

This means that the largest upward spike of height $2\lambda\sqrt{g_2}\log N$ has width of order $N^{1/2-\lambda/2}$.

- The spatial distribution of spikes displays a non-trivial multi-fractal structure. In particular, let $\mathcal{S}_N(\lambda) \stackrel{\text{def}}{=} \left\{ i \in \Lambda_{\varepsilon N} \,:\, \varphi_i \geq 2\lambda\sqrt{g_2}\log N \right\}$ be the set of $\lambda$-*high points*. Then the number of $\lambda$-high points asymptotically satisfies

$$\lim_{N\to\infty} \frac{\log |\mathcal{S}(\lambda)|}{\log N} = 2(1 - \lambda^2) \qquad \text{in probability,}$$

which means that $|\mathcal{S}(\lambda)|$ is of order $N^{2-2\lambda^2}$. In particular, if $\lambda < \gamma < 1$, the average number of $\lambda$-high points in boxes of the form $B_\gamma(i) \stackrel{\text{def}}{=} \{\|j-i\|_\infty \leq N^\gamma\}$, $i \in \Lambda_{\varepsilon N}$, is of order $N^{2\gamma} \cdot N^{2-2\lambda^2}/N^2 = N^{2\gamma - 2\lambda^2}$. However, it is proved that

$$\lim_{N\to\infty} \max_{i \in \Lambda_{\varepsilon N}} \mathrm{P}_{\Lambda_N}\left( \left| \frac{\log |\mathcal{S}(\lambda) \cap B_\gamma(i)|}{\log N} - 2\gamma\left(1 - (\lambda/\gamma)^2\right) \right| > \delta \right) = 0,$$

and thus the number of $\lambda$-high points in a *typical* square box $B_\gamma(i)$ is of order $N^{2\gamma-2\lambda^2/\gamma} \ll N^{2\gamma-2\lambda^2}$. This indicates that the $\lambda$-high points cluster together, as is confirmed by this other estimate,

$$\lim_{N\to\infty} \max_{i\in\Lambda_{\varepsilon N}} P_{\Lambda_N}\left(\left|\frac{\log|\mathcal{S}(\lambda)\cap B_\gamma(i)|}{\log N} - 2\gamma\left(1-\lambda^2\right)\right| > \delta \,\Big|\, i\in\mathcal{S}_N(\lambda)\right) = 0,$$

which shows that the number of $\lambda$-high points in a typical square box $B_\gamma(i)$, *centered on a $\lambda$-high point*, is of order $N^{2\gamma-2\lambda^2\gamma} \gg N^{2\gamma-2\lambda^2}$. Reference [34] also contains interesting informations on the pairwise spatial distribution of high points, but I do not discuss this here.

The geometry and distribution of spikes in dimensions $d \geq 3$ is much simpler. Actually, it turns out that they have the same behavior at leading order as the extremes of the field of i.i.d. Gaussian random variables with the same variance. In particular, their spatial distribution is uniform, and the spikes are thin (*i.e.* their size is logarithmic).

**Pinning by a single site.** A crucial property of low-dimensional ($d=1$ or 2) continuous interfaces is that the local variance of the field has a slow growth. In particular, it turns out that pinning a single point in sufficient to "localize" the field, in the sense that an infinite-volume Gibbs measure exists. More precisely, let us consider the Gaussian measure $P_{\Lambda_N\setminus\{0\}}$, *i.e.* the Gaussian model with 0 boundary conditions outside $\Lambda_N$ *and at the origin*. Then it follows from the random walk representation that, for any $i \in \mathbb{Z}^d$,

$$\lim_{N\to\infty} \mathrm{var}_{P_{\Lambda_N\setminus\{0\}}}(\varphi_i) \asymp \begin{cases} |i| & (d=1), \\ \log|i| & (d=2). \end{cases}$$

Actually, one even has that

$$\sup_{i\neq 0} \frac{\mathrm{var}_{P_{\mathbb{Z}^d\setminus\{0\}}}(\varphi_i)}{\mathrm{var}_{P_{\Lambda(i)}}(\varphi_i)} \asymp 1,$$

where $\Lambda(i) \stackrel{\mathrm{def}}{=} \{j \in \mathbb{Z}^d \,:\, |j-i|_\infty \leq |j|_\infty\}$. Therefore, up to a multiplicative constant, the variance of $\varphi_i$ under $P_{\Lambda_N\setminus\{0\}}$ is the same as it would be in a finite box of radius $|i|_\infty$. (Actually, in two dimensions, it would be better to say that the variance of $\varphi_i$ under $P_{\Lambda_N\setminus\{0\}}$ is the same as it would be in a finite box of radius $|i|_\infty^2$, as the two expressions become then asymptotically equal.)

**Remark 9.** *I have taken 0-boundary conditions outside $\Lambda_N$, but any boundary conditions not growing too fast with $N$ would have given the same result. This shows that pinning of a single point in low dimensions can screen the behavior at infinity.*

**Gradient Gibbs states.** The non-existence of infinite-volume Gibbs states in low dimensions is unfortunate, as it prevents the use of a lot of existing technology. However, we have seen that pinning the field at the origin suffices to have a well-defined thermodynamic limit in any dimension. The resulting field is completely characterized by its gradients (since the value of the field at the origin is known to be zero). These infinite-volume random fields of gradients are called gradient Gibbs measures, or Funaki-Spohn states. They have been studied in details, under assumption of uniform strict convexity, in [56]. Among their results is the fact that for every tilt (*i.e.*, every direction of the mean normal to the interface) there is a unique shift-invariant, ergodic Gibbs measure. The convexity assumption turns out to be necessary, as otherwise it is possible to have coexistence of several phases with same tilt, but different degrees of roughness (*i.e.*, different variances) [10].

The analysis of these gradient Gibbs states provides useful technical tools and has led to several results on the large-scale limit of these fields.

**Convergence to the continuous Gaussian free field.** Let $f$ be a continuous real-valued function with compact support in $D$, a non-empty, compact subset of $\mathbb{R}^d$. Let $\varphi^N$ denote a realization of the Gaussian field under $\mathrm{P}_{ND}$. The action of $\varphi^N$ on $f$ can then be defined as

$$(\varphi^N, f) \stackrel{\text{def}}{=} N^{-d} \sum_{i \in ND} f(i/N) \varphi_i^N .$$

It can then be proved that the sequence $(\varphi^N, f)$ converges in law to the *continuous Gaussian (massless) free field* $\Phi$, which is the centered Gaussian family $(\Phi, f)$ indexed by functions $f$ as above, such that

$$\mathrm{cov}((\Phi, f), (\Phi, g)) = \int_{V \times V} f(x) g(y) \, \mathcal{G}_D(x, y) \, \mathrm{d}x \mathrm{d}y ,$$

where $\mathcal{G}_D$ is the Green function of the Brownian motion killed as it exits $V$.

Actually convergence to a continuous Gaussian free field with suitable covariance also holds for any interaction $V$ which is uniformly strictly convex [81].

A more detailed introduction to the Gaussian free field can be found in [87]. Additionally to its describing the continuum limit of effective interface models, this object has many remarkable properties, including conformal invariance. In particular, it has recently triggered much interest due to its relations to SLE processes.

### *1.3. Discrete effective interface models: basic properties*

Let us now turn to the case of discrete effective interface models. It turns out that they have very different behavior in different situations.

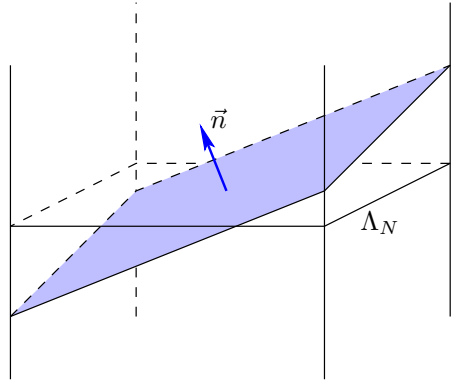

FIG 3. *A sketch of the cylinder $\Lambda_N \times \mathbb{Z}$ and the associated tilted boundary condition $\psi_{\vec{n}}$.*

### 1.3.1. Dimension 1

In dimension 1, continuous and discrete effective interface models should always have identical behavior, at least qualitatively. This is easy to understand, since it essentially reduces to replacing a random walk on the integers by a random walk on the real line, and both have Brownian asymptotics.

### 1.3.2. Dimension 2: roughening transition

The behavior in dimension 2 crucially depends on the boundary conditions and on the temperature. In this review, only the case of tilted boundary conditions will be discussed, *i.e.* boundary conditions $\psi_{\vec{n}}$ approximating a hyperplane through the origin with normal $\vec{n} \in \mathbb{R}^{d+1}$ (not horizontal), see Fig. 3. It is then expected that for all $\vec{n}$ non-vertical, the large-scale behavior of these random interfaces is identical to that of their continuous counterpart. In particular they should have Gaussian asymptotics. This turns out to be quite delicate, and the only rigorous works related to this problem, of which I am aware are [73, 74, 57]. In these papers, the infinite-volume gradient field (see the discussion on gradient Gibbs measures in Subsection 1.2.6), corresponding to the SOS interface with a normal in general non-vertical direction, is analyzed at $\beta = \infty$, and is proved to converge to the continuous Gaussian free field (see Subsection 1.2.6). These results rely on the relationship between these random surfaces and domino tilings.

I am not aware of a single rigorous proof for finite $\beta$, not even of the delocalization of the interface (one might think that at least delocalization should easily follow from these 0-temperature results, since it seems intuitive that thermal fluctuations should increase fluctuations, but there is a delicate problem of interchange of limits, and there are examples in which fluctuations at zero-temperature are much larger than at finite temperatures, see [11]).

The behavior when $\vec{n}$ is vertical, *i.e.* the case $\psi \equiv 0$, is more remarkable.

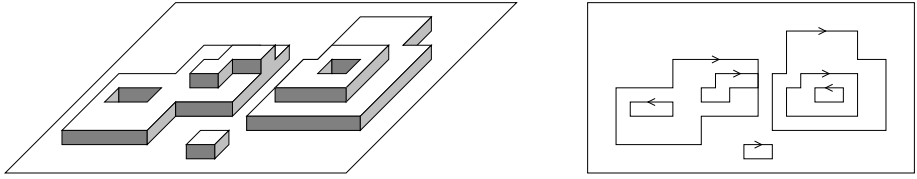

FIG 4. *The oriented level lines of a discrete effective interface; the orientation specify the type of the contour, i.e. whether it is increasing or decreasing the surface height. For convex interactions V, the resulting contours satisfy a Peierls estimate and are therefore amenable to a rigorous study through perturbative arguments.*

It is expected that there is a phase transition, the *roughening transition*, at an inverse temperature $0 < \beta_r < \infty$ such that:

- For all $\beta > \beta_r$, the interface is localized and massive;
- For all $\beta < \beta_r$, the interface has a Gaussian behavior.

The behavior at very large $\beta$ is very well understood. Thanks to the discrete nature of these models, contour techniques are available (see Fig. 4) and most questions can be answered in great details using suitable versions of Pirogov-Sinai theory. Actually, at very large $\beta$, the picture of discrete horizontal interfaces turns out to be the following: The interface is given by the plane $z = 0$ (which is nothing but the ground-state of the model) perturbed by small local fluctuations. In particular, one can easily show that the variance is bounded, the mass positive, the spikes are thin, etc.

The small $\beta$ regime, on the other hand, is still poorly understood. There are no result concerning the Gaussian asymptotics. The stronger results known to date are those given in the celebrated (and difficult!) work of Fröhlich and Spencer [53], who proved that, at small enough $\beta$,

$$\mathrm{var}_{P_\Lambda^\beta}(\varphi_i - \varphi_j) \asymp \log \|j - i\|_2 ,$$

for any $i, j$ and $\Lambda$ large enough. Their result holds for both SOS and DG models.

**Open Problem 3.** *Prove the existence of a roughening transition directly in terms of the interface, without using the mapping to a Coulomb gas, as done in [53]. Prove convergence to a Gaussian field in the rough phase, maybe at least in some limit as $\beta \to 0$.*

*1.3.3. Dimension 3 and higher*

The behavior in dimensions 3 and higher is expected to be radically different: For any $\beta > 0$, the interface should be localized and massive. This was proved in [65] for the DG model in dimension 3 (note that localization is not very surprising since the same is also true for continuous effective interface models; it is the exponential decay of correlations that is remarkable).

*1.4. Additional results*

**Comparison between discrete and continuous models.** It would be interesting to have general inequalities between discrete models and their continuous counterparts. For example, it seems rather plausible that the fluctuations of the discrete models should be dominated by those of the continuous ones. There are indeed a few comparison inequalities of this type, but they are restricted to $V(x) = \frac{1}{2}x^2$ [52] (the subscripts "D" and "C" serve to distinguish between the discrete and continuous models):

$$\text{var}_{P^\beta_{\Lambda,D}}(\varphi_0) \leq \text{var}_{P^\beta_{\Lambda,C}}(\varphi_0), \tag{13}$$

$$E^\beta_{\Lambda,D}(e^{c\varphi_0}) \leq E^\beta_{\Lambda,C}(e^{c\varphi_0}), \tag{14}$$

for any $c \geq 0$. Actually these inequalities also hold in presence of a mass term, see Section 2 for the definition.

**Roughening in one dimension.** It is also possible to study the roughening transition in one-dimensional models. It can be proved that the one-dimensional DG model with $p(i) \sim |i|^{-r}$ describes a rigid interface at any temperatures if $1 < r < 2$, while it is always rough when $r > 2$. The marginal case $r = 2$ has been studied in [54], where it is proved that there is a roughening transition from a rigid to a rough phase as the temperature increases. Moreover in the rough phase, $\text{var}_{P^\beta}(\varphi_i - \varphi_j) \geq c(\beta) \log \|j - i\|_2$.

**Membranes.** Beside the one described above, there is another important class of effective models used to describe *membranes*, whose main representative is still Gaussian, with covariance matrix given by $(\Delta)^{-2}$, instead of $\Delta^{-1}$ as for effective interface models. This change gives rise to radically different properties (especially, these surfaces display huge fluctuations, known in this context as *undulations*, which have dramatic effects on all (de)localization properties). Unfortunately, they are much less tractable from the mathematical point of view, due to the lack of nearly all the main tools used in the study of effective interfaces: no nice random walk representation, no FKG inequalities, etc. I refer to [85, 75] for rigorous results (concerning the entropic repulsion phenomenon) on such models in dimensions 5 and more.

## 2. Massive model

*2.1. Description of the model*

A very crude way to localize the interface is by adding a mass term to the Hamiltonian. This is a very well-known problem, which is particularly easy to analyze in the continuous Gaussian case thanks to the fact that the measure remains Gaussian. I briefly state the results, to allow comparison with those described in other sections. I limit the discussion to the case of the Gaussian model.

Let $\mu \in \mathbb{R}_+$. Let us consider the following perturbation of the measure $P_\Lambda$,

$$P_\Lambda^\mu(d\varphi) \propto \exp\left(-\tfrac{1}{2}\mu^2 \sum_{i\in\Lambda} \varphi_i^2\right) P_\Lambda(d\varphi).$$

The reason for the square in $\mu^2$ in the definition will become clear in a moment.

### 2.2. Main results

The measure $P_\Lambda^\mu$ is still Gaussian, of course, and a similar argument as before yields the following generalization of the random walk representation in this case:

$$\mathrm{cov}_{P_\Lambda^\mu}(\varphi_0, \varphi_i) = \sum_{n\geq 0} \mathbb{P}_i[X_n = j, \tau_\Lambda > n \vee \tau_\mu],$$

where $\tau_\mu$ is a geometric random variable of parameter $\mu^2/(1+\mu^2)$, independent of the random walk. We thus see that, at each step, the random walk has a probability $\mu^2/(1+\mu^2)$ of getting spontaneously killed. This makes it transient in any dimension, and consequently the infinite-volume field is always well-defined. Of course, much more can be said, and precise estimates can be obtained for the variance, mass, maximum, etc. Here I just briefly state some rather rough results, to allow the reader to easily compare with the corresponding ones in other sections: Uniformly in $\mu$ small enough (and under assumption (3) for the last claim)

$$\mathrm{var}^\mu(\varphi_0) \asymp \mu^{-1} \qquad (d=1) \tag{15}$$
$$\mathrm{var}^\mu(\varphi_0) = (\pi\sqrt{\det\mathcal{Q}})^{-1}|\log\mu| + O(\log|\log\mu|) \qquad (d=2) \tag{16}$$
$$P^\mu(\varphi_0 > T) \leq e^{-c(d,\mu)T^2} \qquad (d\geq 1) \tag{17}$$
$$m_{P^\mu}(x) = \mu + o(\mu) \qquad (d\geq 1) \tag{18}$$

The last result shows that $\mu$ is in fact equal (to leading order as $\mu \downarrow 0$) to the mass (defined in (7)). This explains why it was introduced through its square in the definition of $P_\Lambda^\mu$.

Statements of explicit formulas for some of the above quantities can be found in [55, Section 3.3].

## 3. The interface above a wall: entropic repulsion

### 3.1. Description of the model

In this section, I briefly recall what is known about the interaction of an interface with a neutral hard wall, *i.e.* the phenomenon of entropic repulsion; since excellent reviews about this (and related) topic can be found in [13, 60] (for continuous effective models) and [26] (for discrete ones), our discussion will stay rather superficial.

The presence of the hard wall at the sites of $\Lambda_M$ is modeled by the positivity constraint $\Omega_{M,+} \stackrel{\text{def}}{=} \{\varphi_i \geq 0, \forall i \in \Lambda_M\}$. The measure describing this process is then the conditioned measure $P_{\Lambda_N}^{\beta,+(M)}(\,\cdot\,) \stackrel{\text{def}}{=} P_{\Lambda_N}^{\beta}(\,\cdot\,|\Omega_{M,+})$ (I simply write $+$ instead of $+(N)$, when $N=M$).

### 3.2. Main results

When $d=1$, it is very well-known that the interface conditioned to stay above the wall converges to the Brownian excursion. This holds both for discrete and continuous models, and for virtually any interaction $V$ for which the model is well-defined (the corresponding random walk should have increments with bounded variance). Let us therefore turn our attention to what happens in higher dimensions, $d \geq 2$.

I first describe what is known in the case of discrete models. Here one can assume that $N=M$, i.e. that the positivity constraint extends all the way to the boundary of the box. The following estimates are proved in [26] for the SOS and DG models at large values of $\beta$:

$$|\Lambda_N|^{-1} \sum_{i \in \Lambda_N} E_{\Lambda_N}^{\beta,+}(\varphi_i) \asymp (\beta^{-1} \log |\Lambda_N|)^\alpha,$$

where $\alpha = 1$ for the SOS model, and $\alpha = \frac{1}{2}$ for the discrete Gaussian model. The heuristic behind these results is rather simple. We have already seen that these models describe rigid interfaces when $\beta$ is sufficiently large. Assuming that this rigid interface lies at a distance $h$ from the wall, all downward excitations (the downward spikes) of height larger than $h$ are forbidden; this implies that moving the interface from $h$ to $h+1$ provides a gain of entropy of order

$$\prod_{i \in \Lambda_N} \frac{\sum_{k=-h}^{\infty} e^{-\beta k}}{\sum_{k=-h-1}^{\infty} e^{-\beta k}} = \exp\left(|\Lambda_N|O(e^{-\beta h})\right),$$

for the SOS model, and of order $\exp(|\Lambda_N|O(e^{-\beta h^2}))$ for the discrete Gaussian model. On the other hand, there is an associated energetic cost, since one has to raise the interface from height $h$ to height $h+1$; this multiplies the Boltzman weight by $\exp(\beta|\partial \Lambda_N|)$ (for the SOS model) or $\exp(\beta|\partial \Lambda_N|((h+1)^2 - h^2))$ (for the DG model). Balancing these terms yields the claimed result.

**Remark 10.** *To see that this repulsion effect is really due to the downward spikes, it is interesting to compare with what happens in the* wedding cake *model of [3, 4]; see Fig. 5. The latter is a discrete model of random surface which has the following properties: 1. The difference between neighboring heights is 0, 1 or $-1$; 2. a connected region of constant height always has a height larger than that of the region in which it is contained. So, in a sense, fluctuations can only raise the interface, and one might think that it should grow even faster than the SOS or DG models conditioned to be positive. However, a Peierls argument can*

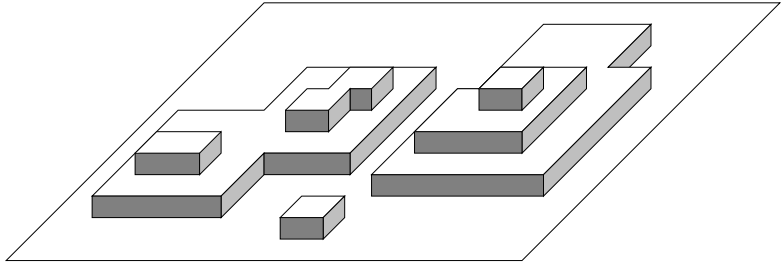

FIG 5. *A picture of the wedding cake model. Gradients can take only values 0 and ±1, and the interface height can only increase. The level lines of this model are exactly identical to the contours of the two-dimensional Ising model.*

*be used to prove that its height is actually finite at low temperature, uniformly in $N$. Notice that in this model, there are* no *downward spikes.*

At higher temperatures, when $d=2$, we have seen that the interface becomes rough, and that it is expected to have Gaussian behavior. In this case, the entropic repulsion effect should be the same as the one for the continuous Gaussian model, to which we turn now, but this has not been proved.

The rest of this subsection is devoted to continuous effective interface models. Let us start with the simpler case of dimensions $d \geq 3$. In this case, it is customary to first take the limit $N \to \infty$, thus studying the measure $\mathrm{P}^{+(M)} \stackrel{\text{def}}{=} \mathrm{P}(\,\cdot\,|\Omega_{M,+})$. One is then interested in the large $M$ asymptotics. It is proved in [17] that, in the case of Gaussian interactions,

$$\lim_{M\to\infty} \sup_{i\in\Lambda_M} \left| \frac{\mathrm{E}^{+(M)}(\varphi_i)}{\sqrt{\log M}} - 2\sqrt{g_d} \right| = 0\,, \qquad (19)$$

where $g_d = \mathrm{var}_{\mathrm{P}}(\varphi_0)$. This shows that the field has a behavior completely similar to its discrete counterpart. Again, it is not to accommodate its typical fluctuations (remember that the variance of the field is finite when $d \geq 3$, and the same can actually be proved for the field conditioned to be positive, by a simple application of Brascamp-Lieb inequality), but rather to accommodate the large downward spikes, exactly in the same way as for the discrete models. Remember (see Subsection 1.2.6) that the typical height of the larger spikes is of order $\sqrt{\log M}$ in a box of size $M$, and therefore of the same order as the repulsion height.

The above estimate remains qualitatively true (that is, without matching upper and lower bounds) for uniformly strictly convex interactions [40].

In the Gaussian case, even more is known about the repelled field: Once the new average is subtracted, it is weakly converging to the unconstrained infinite-volume field [39], which means that both fields look locally the same: There exists a sequence $a_M$, with $\lim_{M\to\infty} a_M/\sqrt{4g_d \log M} = 1$, such that

$$\mathrm{P}^{+(M)}(\,\cdot\, + a_M) \stackrel{M\to\infty}{\Longrightarrow} \mathrm{P}\,.$$

It must be emphasized however that the two fields have the same behavior only *locally*. For example, the minimum of the field in the box $\Lambda_N$ is much smaller than that of the centered repelled field, as can be checked by comparing (19) and (11).

The case $N = M$, i.e. the measure $P_{\Lambda_N}^{+(N)}$ has been treated in [38] and [40]; it turns out that (19) still holds, with the same constant, provided the sup is restricted to $i \in \Lambda_{\epsilon N}$, $0 < \epsilon < 1$. Moreover, estimates for the growth of the interface near to the boundary have also been obtained, showing that the height above a site at distance $L$ from the boundary grows like $\sqrt{\log L}$.

It remains to discuss the two-dimensional case. It is not possible to consider the conditioned infinite-volume measure, since the unconditioned limiting field does not exist. In that case $M$ is chosen such that $M = \epsilon N$ with $0 < \epsilon < 1$. It is then proved in [16] that

$$\lim_{N \to \infty} \sup_{i \in \epsilon \Lambda_N} P_{\Lambda_N}^{+(\epsilon N)} \left( \left| \varphi_i - \sqrt{4g_2} \log N \right| \geq \epsilon \log N \right) = 0, \qquad (20)$$

where $g_2$ is defined in (8). The same, without matching upper and lower bounds, also holds for uniformly strictly convex interactions [40]. In the latter case, the case $M = N$ has also been considered, showing that the same qualitative result also holds in that case; moreover the average height above a site at distance $L$ from the boundary grows like $\log L$.

A few remarks should be made. First of all, it should be noted that the result is qualitatively different from what happens in the two-dimensional DG model, since for the latter the height of the repelled interface has the same order in all dimensions $d \geq 2$. This is due to the fact that the two-dimensional interface is delocalized: As was discussed in Subsection 1.2.6, the spikes in dimension 2 are much fatter than they are in higher dimensions or for the 2-dimensional DG model, and cannot be considered as essentially independent objects; spikes can (and do) actually grow on each other. In order to obtain the sharp result stated above, it is necessary to make a delicate multi-scale analysis. Second, it is interesting to observe that, contrarily to what happens when $d \geq 3$, the repulsion height coincide exactly (at leading order) with the height of the largest downward spikes of the unconstrained field.

A detailed analysis of the downward spikes of the repelled two-dimensional field was done in [34]. It turns out that all the results stated in Subsection 1.2.6 about the spikes of the unconstrained field remain true for the repelled field, both for upward and downward spikes. This indicates that the repelled and original fields are much more alike than in higher dimensions.

**Open Problem 4.** *Prove (or disprove) that the centered repelled field weakly converges to the continuous Gaussian free field.*

**Open Problem 5.** *Prove similar results (even without the proper constants) when $V(x)$ is not uniformly strictly convex. In particular prove it for the continuous SOS model, that is, the continuous effective interface model with interaction* $V(x) = |x|$. Note that the behavior should be different as soon as the tail of $V$ changes. This has important impact on other related issues, as is discussed in Sections 4 and 6.

### 3.3. Additional results

**Pinning by a single site.** Remember that we saw in Subsection 1.2.6 that the field pinned at the origin is always well-defined. One can wonder what happens if one moreover conditions the field to be positive. It turns out that in dimensions 1 and 2 the infinite-volume field still exists, and actually

$$\mathrm{E}^+_{\mathbb{Z}^2\setminus\{0\}}(\varphi_i) \asymp \begin{cases} \sqrt{|i|} & (d=1), \\ \log|i| & (d=2). \end{cases} \quad (21)$$

The situation is however completely different in higher dimensions, since it can be proved that even the expectation of the sites neighboring the origin is infinite! I refer to [60, Lemma 4.4] for nice, probabilistic proofs of this. This is another manifestation of the very different geometry of the spikes in low and high dimensions. Unfortunately, the above result in dimension 2 has only been established for the Gaussian model; an extension to uniformly strictly convex $V$ would also allow an extension of the results of [29], discussed in Section 6, to this class of interactions.

**Open Problem 6.** *Extend* (21) *(when $d = 2$) to the case of uniformly strictly convex interactions $V$.*

**Disordered wall.** In the above, the wall was considered to be perfectly planar. I briefly mention here some studies of this phenomenon in the presence of a rough substrate.

In [7], the wall is modeled by a family of i.i.d. random variables $\sigma_i$, $i \in \mathbb{Z}^d$, $d \geq 3$ independent of the interface $\varphi$. The constraint becomes of course $\varphi_i \geq \sigma_i$, for all $i \in \mathbb{Z}^d$. It is proved that the behavior of the interface depends on the upward-tail of the random variables $\sigma_i$:

- Sub-Gaussian tails: Assume that

$$\lim_{r\to\infty} r^{-2}\mathbb{P}(\sigma_0 \geq r) = -\infty\,.$$

  Then there is no effect and the repulsion height is almost-surely the same as for a flat wall (at leading order), *i.e.* $\sqrt{4g_d \log N}$. This is not so surprising as in this case the extremes of the wall live on a much smaller scale than the repelled field.

- Super-Gaussian tails: Assume there exists $Q > 0$ and $\gamma \in (0, 1)$ such that

$$\lim_{r\to\infty} r^{-2\gamma}\mathbb{P}(\sigma_0 \geq r) = -1/2Q\,.$$

  Then the repulsion height is almost-surely given (at leading order) by $(\sqrt{4Q \log N})^\gamma$, which is much bigger than the repulsion height in the case of a flat wall.

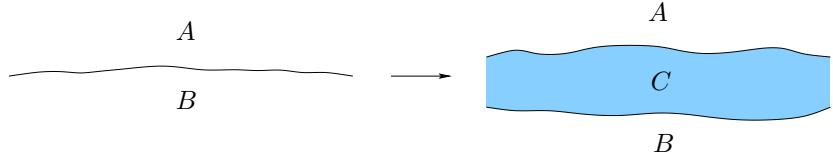

FIG 6. *When $\tau_{AB} > \tau_{AC} + \tau_{BC}$, a layer of the equilibrium phase $C$ is spontaneously created at the interface between the equilibrium phases $A$ and $B$, thus giving rise to a problem of entropic repulsion between multiple interfaces.*

- Almost-Gaussian tails: Assume that there exists $Q > 0$ such that
$$\lim_{r \to \infty} r^{-2} \mathbb{P}(\sigma_0 \geq r) = -1/2Q.$$
Then the repulsion height is almost-surely given (at leading order) by $\sqrt{4(g_d + Q) \log N}$.

In [8], the wall was itself sampled according to the law of a Gaussian effective interface model. It is then proved that, at leading order, this strongly correlated substrate gives rise almost-surely to precisely the same repulsion height as the i.i.d. wall with the corresponding almost-Gaussian upward-tail. This should not be too surprising as only the extreme values of the field modelling the wall play a role, and those have the same behavior (at leading order) in both cases.

**Several interfaces.** Another problem of interest is the analysis of the entropic repulsion effect for several interface (with or without a wall). This models situations in which there are more than two equilibrium phases. For example, one might have three equilibrium phases $A$, $B$ and $C$ and one wishes to analyze the coexistence of phases $A$ and $B$. It may then turn out that it is more favourable for the system to create a layer of phase $C$ between phases $A$ and $B$ (this will happen if the corresponding surface tensions satisfy the inequality $\tau_{AB} > \tau_{AC} + \tau_{BC}$); see Fig. 6.

The case of two Gaussian fields $\varphi^1$ and $\varphi^2$, with the constraint that $\varphi_i^2 \geq \varphi_i^1$, for all $i$, can be mapped to the case of one interface and a flat hard-wall by making a simple change of variables, as long as the covariance matrices of the two fields commute. For example, if both field have a covariance matrix given by the inverse of the $d$-dimensional discrete Laplace operator, then one can consider $\bar\varphi^1 = (\varphi^2 - \varphi^1)/\sqrt{2}$ and $\bar\varphi^2 = (\varphi^1 + \varphi^2)/\sqrt{2}$. Then the constraint reduces to $\bar\varphi^1 \geq 0$, and the two new fields are still Gaussian and independent, so the results of the present section apply directly. I learned from G. Giacomin, however, that the case of two fields with non-commuting covariance matrices, which he has treated [59], is more subtle, and actually requires a new argument.

In [9], the authors consider two $d$-dimensional, $d \geq 3$, Gaussian interfaces $\varphi^1$ and $\varphi^2$, with the constraint that $\varphi_i^2 \geq \varphi_i^1 \geq 0$ for all $i$. The main result is that the height of $\varphi^1$ is still given (at leading order) by $\sqrt{4g_d^1 \log N}$ and is therefore unaffected by the presence of $\varphi^2$ (there is no "pressure" on the lower interface from the upper one). The height of $\varphi^2$ itself is given by $(\sqrt{4g_d^1} + \sqrt{4(g_d^1 + g_d^2)})\sqrt{\log N}$.

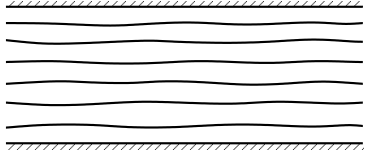

FIG 7. *A stack of interfaces between two fixed walls. Another variant consist in an infinite "one-dimensional" gas of d-dimensional interfaces with positive fugacity.*

In [86], the corresponding question was considered for $K \geq 2$ two-dimensional Gaussian interfaces above a hard wall. The results obtained are completely similar.

## 4. Confinement between two walls

### *4.1. Description of the model*

The case of several interfaces considered at the end of the previous section has the property that the interfaces have all the room they wish to move away from each other. There are important situations where this is not the case (for example, when modelling commensurate/incommensurate transitions, see [26] for more on this). In that case one would like to study a stack of interfaces with fixed separation (in other word, one would like to study a "one-dimensional" gas of $d$-dimensional interfaces with fixed density, or fugacity); see Fig. 7. This is a much more difficult problem, and essentially only a very simple caricature of this situation has been studied rigorously so far: the case of a single interface between two walls. It turn out that this is a nice model (in the Gaussian case, anyway), as its critical behavior can be analyzed in enough details to prove deviations from mean-field in low dimensions; see Section 5 on pinning for another one of the few situations where this can be done.

Let $\ell > 0$. Let us consider the following modification of our basic Gibbs measure
$$\mathrm{P}_{\Lambda;\ell}^{\beta}(\,\cdot\,) \stackrel{\text{def}}{=} \mathrm{P}_{\Lambda}^{\beta}(\,\cdot\,|\,|\varphi_i| \leq \ell, \forall i \in \mathbb{Z}^d)\,.$$
This models an interface constrained to lie inside a horizontal slab of height $2\ell$.

### *4.2. Main results*

**Continuous effective interface models.** Of course, once constrained inside a slab, the interface is localized in any dimension. It has also been proved in [79] that it becomes massive in any dimension, but with an estimate for the mass that is only correct in dimensions 3 and above. This result has later been improved in order to get the qualitatively correct behavior for the mass and related quantities in [26]. Their results, which are restricted to the nearest-neighbor Gaussian model, can be stated as follows: Set $\nu(2) = 1$ and $\nu(d) = 2$ if $d \geq 3$, then [2]

- The mass of the confined interface has the following large-$\ell$ asymptotics:

$$m_{\mathrm{P}_\ell} = \begin{cases} O(\ell^{-2}) & (d=1), \\ \exp(-O(\ell)) & (d=2), \\ \exp(-O(\ell^2)) & (d \geq 3). \end{cases} \quad (22)$$

- The variance of the confined interface has the following large-$\ell$ asymptotics:

$$\mathrm{var}_{\mathrm{P}_\ell}(\varphi_0) \asymp \begin{cases} \ell^2 & (d=1), \\ \ell & (d=2), \end{cases} \quad (23)$$

and satisfies $0 \leq \mathrm{var}_{\mathrm{P}}(\varphi_0) - \mathrm{var}_{\mathrm{P}_\ell}(\varphi_0) \leq \exp(-O(\ell^2))$, for $d \geq 3$.
- The probability that the Gaussian interface remains inside the slab of height $2\ell$ has the following asymptotic large-$\ell$ behavior,

$$|\Lambda_N|^{-1} \log \mathrm{P}_{\Lambda_N}(|\varphi_i| \leq \ell, \forall i \in \Lambda_N) = \begin{cases} O(\ell^{-2}) & (d=1), \\ \exp(-O(\ell)) & (d=2), \\ \exp(-O(\ell^2)) & (d \geq 3), \end{cases} \quad (24)$$

for all $N > N_0(\ell)$.

The last result has very recently been given a sharper form in dimensions $d \geq 3$ in [84]; it turns out, unsurprisingly, that this probability has the same leading-order behavior as the corresponding i.i.d. Gaussian field,

$$|\Lambda_N|^{-1} \log \mathrm{P}_{\Lambda_N}(|\varphi_i| \leq \ell, \forall i \in \Lambda_N) = \exp(-\frac{\ell^2}{2g_d}(1+o(1))).$$

The lower bound is actually an immediate consequence of Griffiths' inequality.

**Remark 11.** *The estimate (24) can be interpreted as an estimate of the asymptotic behavior of the free energy of the constrained interface.*

**Remark 12.** *A completely different proof of (24), valid for arbitrary uniformly strictly convex interactions, can be found in [90]; it makes use of the results of the preceding section on entropic repulsion, together with correlation inequalities. However, it does not permit to recover the estimates (23) and (22) for the variance and the mass of the confined interface.*

---

[2]The notation $O(b)$ denotes a function such that $c_- b \leq O(b) \leq c_+ b$, for all large (or small depending on the context) $b$, for some constants $0 < c_- \leq c_+ < \infty$, possibly depending on fixed parameters (dimension, temperature, etc.).

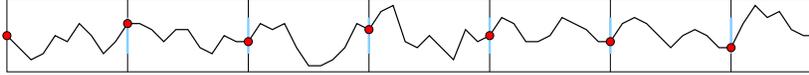

FIG 8. *The probability that a random walk stays in a slab of height $2\ell$ for a time $N$ can be estimated from below as follows. One splits the time-interval in pieces of length $\ell^2$ and force the walk to come back to the interval $[-\ell/2, \ell/2]$ at the end of each piece; this has a probabilistic cost $\exp(-O(N/\ell^2))$. Inside each piece the desired event has strictly positive probability, uniformly in $\ell$, which gives an additional factor $\exp(-O(N/\ell^2))$. An analogous reasoning yields the corresponding upper bound; see, e.g., [68].*

**Open Problem 7.** *Prove the claims* (23) *and* (22) *on the variance and the mass beyond the Gaussian case. Also, in the Gaussian case, try to get sharper results.*

Of course, all these three results should be closely connected. Indeed, once the mass has been computed, one should be able to replace the slab-constraint by this effective mass, and all computations become purely Gaussian. Doing this does indeed yield the above estimates for the variance (just plug the above expression for the mass into (15) or (16)) and the free energy. This is actually the way the proof in [26] works.

It is simple to understand heuristically the scaling (24) of the free energy when $d = 1$, see Fig. 8; actually this can easily be turned into a rigorous proof. To understand (24) in higher dimensions, observe that the maximum of the Gaussian field over the box $\Lambda_N$ is of order $\log N$ ($d = 2$), resp. $\sqrt{\log N}$ ($d \geq 3$), and therefore the typical distance over which the interface feels the confinement should be given by $\log N \sim \ell$ ($d = 2$), resp. $\log N \sim \ell^2$ ($d \geq 3$). So one essentially has to pay a fixed price for each portion of the box of linear size $\exp(\ell)$ ($d = 2$), resp. $\exp(\ell^2)$ ($d \geq 3$), which is precisely (24). The proof in [90] goes along these lines.

**Discrete effective interface models.** The corresponding results for discrete heights have also been established in [26]. For $d = 1$, and any interactions $V$, one has the same results as before. So only the higher-dimensional ones are discussed here. The most detailed results concern the free energy in the case of the DG model:

$$|\Lambda_N|^{-1} \log P^\beta_{\Lambda_N}(|\varphi_i| \leq \ell, \forall i \in \Lambda_N) \begin{cases} = \exp(-O(\ell^2)) & (d = 2, \beta \gg 1), \\ \leq \exp(-O(\ell)) & (d = 2, \forall \beta), \\ = \exp(-O(\ell^2)) & (d \geq 3, \forall \beta). \end{cases}$$

In the case of the SOS model, it is only proved that when $d = 2$ and $\beta \gg 1$, the above quantity is of order $\exp(-O(\ell))$. This is one example of the fact that in phenomena depending on the behavior of spikes, the tail behaviour of the interaction is crucial.

There are only few results about the mass for the DG model. Basically, it is only known that, unsurprisingly, the mass is positive for $\ell < \infty$ if it is already positive for $\ell = \infty$ (remember that this is the case when $d = 2$ and $\beta$ is large, or $d \geq 3$ and $\beta$ is arbitrary). The only result pertaining to the most interesting case of $d = 2$ and $\beta$ small is that the covariances are summable (*i.e.*, the susceptibility is finite), for any $d \geq 1$, $\beta < \infty$ and $\ell < \infty$.

### *4.3. Additional results*

**Centering of the interface in the slab.** In this section, I have always considered 0-boundary conditions. However, it can be proved (see [26]) that for the DG and SOS models at large $\beta$, as well as for the continuous Gaussian model in $d \geq 3$, all boundary conditions lead to the same infinite-volume Gibbs state, as long as $\ell < \infty$. This is due to a "centering" of the interface in the middle of the slab, which can heuristically be understood in a way similar to what we saw for entropic repulsion: a simple computation shows that there is loss of entropy, due to the spikes, when the interface is not centered.

**Two confined interfaces.** The case of two $d$-dimensional Gaussian interfaces, $d \geq 3$, has recently been studied [84]. The result is that, when the two interfaces have the same covariance matrix, the ratio between the average distance between the two interfaces and between each of them and the closest wall is $\sqrt{2}$ at leading order.

## 5. Pinning

### *5.1. Description of the model*

In this section, we are going to investigate the effect of a self-potential favouring a finite neighborhood of zero. Namely, letting $a > 0$, $b > 0$, and $W_{a,b}(x) = -b\,\mathbf{1}_{\{|x| \leq a\}}$, let us consider the following perturbed probability measure:

$$\mathrm{P}_\Lambda^{a,b}(\mathrm{d}\varphi) \propto \exp\Big(-\sum_{i \in \Lambda} W_{a,b}(\varphi_i)\Big)\,\mathrm{P}_\Lambda(\mathrm{d}\varphi)\,, \tag{25}$$

and similarly in the case of discrete heights.

There are several reasons why this an interesting problem: its relation to the wetting phenomenon, which can be seen as a competition between localization by such a self-potential and entropic repulsion, see Section 6; the fact that, in the case of continuous effective interface models, it provides a weak perturbation of the measure $\mathrm{P}_\Lambda$ breaking the continuous symmetry ultimately responsible for the unbounded fluctuations of the interface in low dimensions, as explained in Subsection 1.2.3; and the fact that the corresponding non-trivial critical behavior can be rigorously studied in great details.

Historically, its analysis for one-dimensional effective interface models was initiated by a desire to better understand the wetting phenomenon discovered through exact computations by D. Abraham in the two-dimensional Ising model [1][3]; more on that in Section 6. It was thus considered useful to prove a similar result in the simpler context of one-dimensional effective interface models. Of particular interest to these earlier works [31, 89, 28, 2] was the dramatic difference in behavior between cases where the localizing self-potential was coupled or not with a positivity constraint: in the former case there is a delocalized phase for weak enough self-potential, while localization always occur in the latter case. After these early works limited to the one-dimensional model, the first rigorous work I am aware is [27] in which the authors introduced a new random walk representation (different from the one described in the introduction) and applied it in particular to prove exponential decay of the 2-point function under the measure $P^{a,b}$ in dimensions $d \geq 3$ for the Gaussian model. A study of the much more delicate two-dimensional case was then done by Dunlop *et al.* [48], who proved that the field is localized for any strictly positive values of $a$ and $b$, in the sense that $E^{a,b}(|\varphi_0|) < \infty$. This result was improved to show finiteness of the variance, as well as exponential decay of the 2-point function by Bolthausen and Brydges [15]. All these results were limited to Gaussian interactions. A more general approach was then developed in [43] in order to treat the case of non-Gaussian (but uniformly strictly convex) interactions and, as a side-product, provided stronger results such as the correct tail for the one-site marginals. This method was then successively improved, first in [70] to prove exponential decay of the 2-point function for this class of models, and then in [20] to establish precise estimates on the critical behavior in any dimensions for possibly long-range Gaussian interactions.

Concerning discrete effective interface models, the situation is as follows: in one-dimension, there is no difference with the case of continuous heights, which is discussed below; in dimensions 3 and more, since the interface is expected (and, for the DG model, proved) to be localized and massive at all temperatures, the pinning potential has no interesting effect. The only interesting case is thus the two-dimensional one. Of course, at low temperatures, the interface is also localized and massive, so the pinning potential is irrelevant At higher temperature, however, we know that the free model should undergo a roughening transition (as proved for the SOS and DG models). One would then expect that the presence of an arbitrarily weak pinning potential will localize and render massive the corresponding interface. However, as we'll see below, the interface gets more and more weakly localized as the pinning strength goes to zero (the variance diverges, the mass vanishes), so one might still expect a transition from a microscopically flat interface to a weakly localized one, although it is not clear whether this should happen smoothly, or through a real phase transition. There are no results on this model unfortunately.

**Open Problem 8.** *What happens to the roughening transition of the two-dimensional DG model in the presence of a pinning potential?*

---

[3]Actually, this was done earlier by McCoy and Wu [80], but they failed to interpret properly what they had computed.

In the rest of this section, we restrict our attention to the (continuous) Gaussian measure, with a few remarks on what is known in the more general case of uniformly strictly convex interactions. Although most of what is discussed below also holds for $W_{a,b}$, I only discuss a particular limiting case of the measure $P_\Lambda^{a,b}$, having the advantage of being nicer from a mathematical point of view, and of depending only on a single parameter. Let $\eta > 0$; the measure with $\delta$-*pinning* is defined by the following weak-limit

$$P_\Lambda^\eta \stackrel{\text{def}}{=} \lim_{\substack{a\to 0, b\to\infty \\ 2a(e^b-1)\to\eta}} P_\Lambda^{a,b}.$$

The latter can also be written more explicitly as

$$P_\Lambda^\eta(\mathrm{d}\varphi) = (Z_\Lambda^\eta)^{-1} \exp\left(-H_\Lambda(\varphi)\right) \prod_{i\in\Lambda} (\mathrm{d}\varphi_i + \eta\delta_0(\mathrm{d}\varphi_i)) \prod_{j\notin\Lambda} \delta_0(\mathrm{d}\varphi_j), \qquad (26)$$

where the partition function $Z_\Lambda^\eta$ is the normalization constant, $\delta_0$ is the Dirac mass at 0, and $H_\Lambda(\varphi)$ is the Hamiltonian (1).

### 5.2. Mapping to a problem of RWRE

To analyze the properties of this model, it is very useful to first map it onto a model of random walk in an annealed random environment of killing obstacles. This is done by simply expanding the product in (26) and using the random walk representation (5) [4].

First observe that for any function $f$

$$E_\Lambda^\eta(f) = (Z_\Lambda^\eta)^{-1} \int f(\varphi) \exp\left(-H_\Lambda(\varphi)\right) \prod_{i\in\Lambda} (\mathrm{d}\varphi_i + \eta\delta_0(\mathrm{d}\varphi_i))$$
$$= \sum_{A\subset\Lambda} \eta^{|A|} \frac{Z_{\Lambda\setminus A}}{Z_\Lambda^\eta} E_{\Lambda\setminus A}(f),$$

and thus the measure $P_\Lambda^\eta$ is nothing but a convex combination of measures $P_{\Lambda\setminus A}$ indexed by subsets of the box $\Lambda$ weighted according to the probability measure $\zeta_\Lambda^\eta(A) \stackrel{\text{def}}{=} \eta^{|A|} Z_{\Lambda\setminus A}/Z_\Lambda^\eta$. Now, combining this with the random walk representation (5), we obtain

$$E_\Lambda^\eta(\varphi_i\varphi_j) = \sum_{A\subset\Lambda} \zeta_\Lambda^\eta(A) \sum_{n\geq 0} \mathbb{P}_i[X_n = j,\, \tau_{\Lambda\setminus A} > n], \qquad (27)$$

where I recall that $\tau_B = \min\{n : X_n \notin B\}$. We thus see that the control of the covariance of the field is reduced to the analysis of the Green function of the random walk $X$ in an annealed random environment of killing traps distributed according to $\zeta_\Lambda^\eta$. It has become standard terminology to call these traps *pinned* sites (for obvious reasons). Let us denote by $\mathcal{A}$ the random subset of $\Lambda$ distributed according to $\zeta_\Lambda^\eta$.

---

[4] In the case of the square-well potential $W_{a,b}$, the expansion of the product in (26) is replaced by the expansion of the product of terms of the form $e^{b\mathbf{1}_{\{|\varphi_i|\leq a\}}} = (1 + (e^b - 1)\mathbf{1}_{\{|\varphi_i|\leq a\}})$.

Of course, the first step in order to turn the above representation into something useful is to obtain informations about the behavior of the probability measure $\zeta_\Lambda^\eta$.

### 5.3. The distribution of pinned sites

In this subsection I briefly describe what is known about the distribution $\zeta_\Lambda^\eta$ of the traps. Let us denote by $\nu_\Lambda^\rho$ the restriction of the Bernoulli site percolation process with density $\rho$ to the box $\Lambda$. Then it is proved in [20] that there exist constants $0 < c_1(d) < c_2(d) < \infty$ such that, for any $\Lambda$, any $B \subset \Lambda$, and any $\eta > 0$ small enough,

$$\nu_\Lambda^{\rho_+(d)}(\mathcal{A} \cap B = \emptyset) \leq \zeta_\Lambda^\eta(\mathcal{A} \cap B = \emptyset) \leq \nu_\Lambda^{\rho_-(d)}(\mathcal{A} \cap B = \emptyset), \qquad (28)$$

where $\rho_\pm(d) \stackrel{\text{def}}{=} \begin{cases} c_\pm(1)\,\eta^2 & (d=1) \\ c_\pm(2)\,\eta |\log \eta|^{-1/2} & (d=2) \\ c_\pm(d)\,\eta & (d \geq 3) \end{cases}$.

Informally, as long as we are only interested in controlling the covariances of the field, this shows that the random environment can be thought of as being Bernoulli. Indeed, (27) can be rewritten as

$$\mathrm{E}_\Lambda^\eta(\varphi_i \varphi_j) = \sum_{n \geq 0} \mathbb{E}_i \left[ \mathbf{1}_{\{X_n = j\}} \mathbf{1}_{\{X_{[0,n]} \subset \Lambda\}} \zeta_\Lambda^\eta(\mathcal{A} \cap X_{[0,n]} = \emptyset) \right],$$

where $X_{[0,n]} = \{X_k : 0 \leq k \leq n\}$ is the range of the random walk. The previous result thus allows us to substitute in the last equation the measure $\zeta_\Lambda^\eta$ by a Bernoulli measure of suitable density.

**Remark 13.** *One might wonder whether it is possible to compare the measures $\zeta_\Lambda^\eta$ and the corresponding Bernoulli in a stronger sense than above. This is indeed possible, but only when $d \geq 3$. If $\mu$ and $\nu$ are two measures on the set of subsets of $\{0,1\}^\Lambda$, for some finite $\Lambda$, let us say that $\mu$ stochastically dominates $\nu$ if $\mu(f) \geq \nu(f)$, for all increasing functions $f$, and let us say that $\mu$ strongly stochastically dominates $\nu$ if $\mu(x \in \mathcal{A} | \mathcal{A} \setminus \{x\} = C) \geq \nu(x \in \mathcal{A} | \mathcal{A} \setminus \{x\} = C)$, for any increasing function $f$ and any $C \subset \Lambda \setminus \{x\}$. Clearly strong stochastic domination implies stochastic domination. When $d \geq 3$, it can be proved that $\zeta_\Lambda^\eta$ is strongly stochastically dominated by $\nu_\Lambda^{\rho_+(d)}$ and strongly stochastically dominates $\nu_\Lambda^{\rho_-(d)}$, see below. This is not true when $d \leq 2$, as can be seen by taking $\Lambda$ a big square centered at the origin and $C = \emptyset$, since in that case the variance of the field at the origin diverges with the size of $\Lambda$. In fact, it can be shown that even simple stochastic domination must fail: There cannot be domination by a Bernoulli process of density $o(\eta)$, see [20] for more details.*

**Remark 14.** *In the case of uniformly strictly convex interactions, both the representation* (27) *and the above results on* $\zeta_\Lambda^\eta$ *hold true[5]. However, in that case the law of the random walk depends on the distribution of pinned points. It is therefore much more delicate to make use of* (28) *than in the Gaussian case. This is the main difficulty in proving the exponential decay of correlations in* [70].

Since these estimates on the distribution of pinned sites constitute the core of the analysis of pinning, I now give some ideas of the proof (in the Gaussian settings, to make things easier); the complete argument can be found in [20].

In dimensions $d \geq 3$, it is actually possible to give the complete proof as it is very elementary. Let $C \subset \Lambda \setminus \{i\}$. Then

$$\zeta_\Lambda^\eta \left(\mathcal{A} \ni i \,|\, \mathcal{A} = C \text{ off } i\right) = \left(1 + \frac{1}{\eta} \frac{Z_{\Lambda \setminus C}}{Z_{\Lambda \setminus (C \cup \{i\})}}\right)^{-1} = \frac{1 + o(1)}{\sqrt{2\pi \mathrm{var}_{P_{\Lambda \setminus C}}(\varphi_i)}} \eta,$$

observing that $Z_{\Lambda \setminus (C \cup \{i\})}/Z_{\Lambda \setminus C}$ is nothing else but the density at 0 of the Gaussian random variable $\varphi_i$ under $P_{\Lambda \setminus C}$. Now, since $d \geq 3$, we know that the variance in the last expression is bounded away from 0 and $\infty$, uniformly in $\Lambda$ and $C$. This immediately implies (28) in that case (actually, this implies the much stronger claim in Remark 13).

Of course, such an argument cannot be used in lower dimensions, since the original field is delocalized, and therefore the variance in the last expression diverges with $\Lambda$. I now describe how to deal with this problem when $d = 2$, in order to prove the upper bound in (28) (the lower bound is much easier). The idea is to prove first (28) when $B$ is composed of cells of a square grid of spacing $K\eta^{-1/2}|\log \eta|^{1/4}$, for some $K$ to be chosen later. I assume that $B$ is connected, to simplify the argument.

Let us write $B^0 \stackrel{\text{def}}{=} B$, and define $B^k$, $k \geq 1$, as being the set of cells obtained by adding to $B^{k-1}$ all its neighboring cells. Let also $\bar{k}$ be the largest value of $k$ such that $B^k \subset \Lambda$, and set $\bar{\mathcal{A}} \stackrel{\text{def}}{=} \mathcal{A} \cup (\mathbb{Z}^2 \setminus \Lambda)$. We can then write

$$\zeta_\Lambda^\eta \left(\bar{\mathcal{A}} \cap B = \emptyset\right) \leq \sum_{k=0}^{\bar{k}} \zeta_\Lambda^\eta \left(\bar{\mathcal{A}} \cap B^k = \emptyset \,|\, \bar{\mathcal{A}} \cap B^{k+1} \neq \emptyset\right).$$

It is therefore sufficient to prove our claim for the conditional probabilities in the RHS, as this implies the same claim with a different constant for the LHS. This shows that there is no loss of generality in assuming that there is at least one point of $\mathcal{A}$ among all cells neighboring $B$; let us denote by $\mathcal{E}$ the corresponding

---

[5] The only estimate where the Gaussian assumption was used in the proof in [20], see Footnote 2 therein, can be extended using Brascamp-Lieb inequality, as shown in [61].

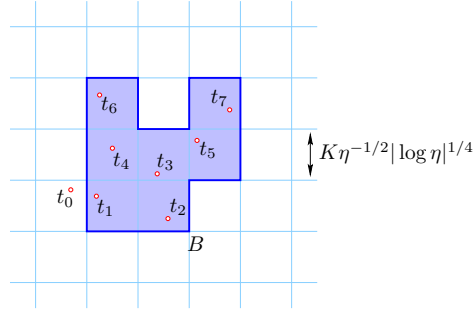

FIG 9. *The set $B$, composed of large cells, with one site $t_k$ in each of its cells. We can assume, without loss of generality, that there is at least one site $t_0$ of $\mathcal{A}$ in a cell neighboring $B$.*

event. We then write

$$\zeta_\Lambda^\eta \left( \mathcal{A} \cap B = \emptyset \,|\, \mathcal{E} \right) = \frac{\sum_{\substack{A \in \mathcal{E} \\ A \cap B = \emptyset}} \eta^{|A|} Z_{\Lambda \setminus A}}{\sum_{A \in \mathcal{E}} \eta^{|A|} Z_{\Lambda \setminus A}}$$

$$= \frac{\sum_{\substack{A \in \mathcal{E} \\ A \cap B = \emptyset}} \eta^{|A|} Z_{\Lambda \setminus A}}{\sum_{C \subset B} \eta^{|C|} \sum_{\substack{A \in \mathcal{E} \\ A \cap B = \emptyset}} \eta^{|A|} Z_{\Lambda \setminus (A \cup C)}}$$

$$\leq \left\{ \sum_{C \subset B} \eta^{|C|} \min_{\substack{A \in \mathcal{E} \\ A \cap B = \emptyset}} \frac{Z_{\Lambda \setminus (A \cup C)}}{Z_{\Lambda \setminus A}} \right\}^{-1}$$

$$\leq \left\{ \sum_{C \in \mathcal{D}(B)} \eta^{|C|} \min_{\substack{A \in \mathcal{E} \\ A \cap B = \emptyset}} \frac{Z_{\Lambda \setminus (A \cup C)}}{Z_{\Lambda \setminus A}} \right\}^{-1}.$$

In the last expression, the sum has been restricted (see Fig. 9) to all

$$C \in \mathcal{D}(B) \stackrel{\text{def}}{=} \{ C \subset B : C \text{ contains exactly one point in each cell of } B \}.$$

We now need to bound from below the ratio of partition functions. This is done by relying on the fundamental property of the two-dimensional interfaces discussed in Subsection 1.2.6: Pinning the interface at one site is enough to bring it down, and in particular to reduce drastically its variance in the neighborhood of this site. To make use of this fact, it is convenient to write the ratio as a telescopic product, so that the problem is reduced to that of the ratio between two partition functions differing by a single pinned site. Let us number the sites of $C = \{t_1, \ldots, t_{|C|}\}$ in such a way that for any $k$, at least one of the sites of $A \cup \{t_1, \ldots, t_{k-1}\}$ belongs to a cell neighboring the one containing $t_k$; this is always possible because $A \in \mathcal{E}$. Using the notation $A_k \stackrel{\text{def}}{=} A \cup \{t_1, \ldots, t_k\}$, $k \geq 1$ and $A_0 \stackrel{\text{def}}{=} A$, the ratio can be rewritten in the following way,

$$\frac{Z_{\Lambda\setminus(A\cup C)}}{Z_{\Lambda\setminus A}} = \prod_{k=1}^{|C|} \frac{Z_{\Lambda\setminus A_k}}{Z_{\Lambda\setminus A_{k-1}}}. \tag{29}$$

Observing, as before, that $Z_{\Lambda\setminus A_k}/Z_{\Lambda\setminus A_{k-1}}$ is the density at 0 of the Gaussian random variable $\varphi_{t_k}$ under the measure $P_{\Lambda\setminus A_{k-1}}$, we have that

$$\frac{Z_{\Lambda\setminus A_k}}{Z_{\Lambda\setminus A_{k-1}}} = \frac{1}{\sqrt{2\pi \text{var}_{P_{\Lambda\setminus A_{k-1}}}(\varphi_{t_k})}}. \tag{30}$$

To estimate the variance, remember that the points of $C$ were numbered in such a way that, under $P_{\Lambda\setminus A_{k-1}}$, there is already a pinned point in a cell neighboring the one containing $t_k$, and thus, by (21), we have

$$\text{var}_{P_{\Lambda\setminus A_{k-1}}}(\varphi_{t_k}) \leq c |\log \eta|.$$

Putting all this together, we obtain

$$\zeta_\Lambda^\eta(\mathcal{A} \cap B) \leq \left\{ \sum_{C \in \mathcal{D}(B)} \eta^{|C|} (c\sqrt{|\log \eta|})^{-|C|} \right\}^{-1}.$$

Since, for all $C \in \mathcal{D}(B)$, $|C| = |B|/(K\eta^{-1}|\log \eta|^{1/2})$, and, moreover,

$$\#\mathcal{D}(B) = (K\eta^{-1}|\log \eta|^{1/2})^{|B|/(K\eta^{-1}|\log \eta|^{1/2})},$$

we see that the claim follows once we choose $K$ large enough.

To prove (28) for arbitrary subsets $B$ is more delicate. The idea is to use the previous result to show that "most" points of $B$ must typically be at a distance at most $O(\eta^{-1}|\log \eta|^{1/2})$ from pinned sites, and therefore the variance of the field cannot be too large and arguments similar to what we did before apply again.

The proof in dimension 1 can be done in the same way.

### 5.4. Main results

Now that we have informations on the distribution of pinned sites, we can turn to the behavior of the field itself. As mentioned earlier, it turns out that an arbitrarily weak $\delta$-pinning is sufficient to localize the interface, in a very strong sense: It is proved in [43] that, for any $\eta$ small enough, and uniformly in all $T$ large enough,

$$-\log P^\eta(\varphi_0 > T) \asymp \begin{cases} T & (d=1) \\ T^2/\log T & (d=2) \\ T^2 & (d \geq 3) \end{cases}$$

In particular the infinite-volume Gibbs measure $P^\eta$ is well-defined in any dimension $d \geq 1$. The heuristic behind this result is elementary (in dimension 2; dimension 1 is completely similar): The probability of having no pinned sites at distance smaller than $R$ from $0$ is $e^{-O(\eta/\sqrt{|\log \eta|})R^2}$, while when this happens $\varphi_0$ is a Gaussian variable with variance of order $\log R$. Combining these two estimates yields the claim above.

In addition to localizing the field in low dimensions, the pinning potential generates a mass in any dimension [70]: For any $d \geq 1$, for any $\eta$ small enough,

$$\inf_{x \in \mathbb{S}^{d-1}} m_{P^\eta}(x) > 0\,.$$

Notice that the tail in dimensions 1 and 2 are *not* Gaussian, which shows that, even though there is localization and exponential decay of correlations, the resulting field is very different from the massive Gaussian field of Section 2.

**Remark 15.** *Although stated here for Gaussian interactions, these results hold for arbitrary uniformly strictly convex interactions.*

Of course, as the intensity $\eta$ of the pinning goes to zero, the localization of the interface becomes weaker and weaker. In order to understand how this progressive delocalization occurs (or, as the physicists say, to analyse the corresponding critical behavior), it is useful to study the variance and the mass of the field. The most precise results in this context are those of [20] and can be stated as follows:

- Consider the one-dimensional model with an arbitrary nearest-neighbor interaction $V : \mathbb{R} \to \mathbb{R}^+$ such that $\int e^{-\frac{1}{2}V(x)}\mathrm{d}x < \infty$, $\int x e^{-\frac{1}{2}V(x)}\mathrm{d}x = 0$, and $\int x^2 e^{-\frac{1}{2}V(x)}\mathrm{d}x / \int e^{-\frac{1}{2}V(x)}\mathrm{d}x = 1$. Then

  $$E^\eta(\varphi_0^2) = \tfrac{1}{2}\eta^{-2} + o(\eta^{-2})\,.$$

  and

  $$m_{P^\eta} = \tfrac{1}{2}\eta^2 + o(\eta^2)\,.$$

- (Gaussian only) Assume that $d = 2$ and (2) holds. Then, for all $\eta$ small enough,

  $$E^\eta(\varphi_0^2) = \left(2\pi\sqrt{\det \mathcal{Q}}\right)^{-1} |\log \eta| + O(\log|\log \eta|)\,. \tag{31}$$

- (Gaussian only) Assume (3) holds. Then, for all $\eta$ small enough,

  $$m_{P^\eta} \asymp \begin{cases} \sqrt{\eta}\,|\log \eta|^{-3/4} & (d = 2) \\ \sqrt{\eta} & (d \geq 3) \end{cases}$$

**Remark 16.** *In dimension 1, only the nearest neighbor case was considered in [20], in order to benefit from the simple renewal structure available in that case. Extensions to long-range interactions should be possible though, since the main result on the trap distribution, see Subsection 5.3, also holds in that case.*

**Remark 17.** *The assumptions on the range of the interaction, in dimensions $d \geq 2$, are essentially optimal. Indeed, when $p(\,\cdot\,)$ has no second moment, the variance of the field may remain finite even when $\eta = 0$ (the random walk $X$ may become transient). Similarly, when $p(\,\cdot\,)$ does not have exponential moments, then the mass is zero for every values of $\eta$, since in the random walk representation of $\mathrm{E}^\eta(\varphi_0 \varphi_i)$, the random walk can go from $0$ to $i$ in a single step, a strategy which obviously has a subexponential cost in this case.*

**Remark 18.** *Again, if one takes the asymptotic behaviour of the mass given above for granted, then one can (on a heuristic level) easily recover the claim on the variance by replacing the model with pinning potential by the Gaussian model with the corresponding effective mass. In particular, injecting $\mu = O(\sqrt{\eta}|\log \eta|^{-3/4})$ into (16), one immediately gets (31), with the correct constant. Of course, this is not a rigorous argument, and the measure with pinning* is *not at all a Gaussian measure, as we already saw at the beginning of this subsection.*

**Remark 19.** *In dimension $d = 2$, in the case of non-Gaussian, but uniformly strictly convex interactions, a corresponding qualitative result for the variance is also known [43]: For $\eta$ small enough,*

$$\mathrm{E}^\eta(\varphi_0^2) \asymp |\log \eta|\,.$$

*Quantitative results in that case, however, would require a much better understanding of the corresponding random walk representation.*

I briefly give some partial proofs for the results when $d \geq 2$. The schemes of proof of these estimates are not very difficult, but it turns out that, in order to get results as sharp as those stated above, rather precise informations on the random walk are necessary. In particular, one needs delicate estimates on its range, going beyond Donsker-Varadhan's asymptotics. As this is not specific to random surfaces and quite technical, this will not be discussed here, and we refer to the appendices in [20] for the details.

**The variance.** I assume that $d = 2$ and only prove the simpler lower bound. This is very easy. Let $B_\eta \stackrel{\text{def}}{=} \{i \in \mathbb{Z}^2 \,:\, \|i\|_\infty \leq \frac{1}{2}\eta^{-1/2}|\log \eta|^{-1/4}\}$. The Green function can be estimated from below by imposing that no trap enter the box $B_\eta$; notice that, according to (28), this box is small enough for this event to have probability close to 1. Assuming that the box $\Lambda$ is very large (we want to take the thermodynamic limit in the end, anyway), we can write

$$\mathrm{E}^\eta_\Lambda(\varphi_0^2) = \sum_{A \subset \Lambda} \zeta^\eta_\Lambda(A) \, \mathrm{E}_{\Lambda \setminus A}(\varphi_0^2)$$

$$\geq \zeta^\eta_\Lambda(\mathcal{A} \cap B_\eta = \emptyset) \min_{A \cap B_\eta = \emptyset} \mathrm{E}_{\Lambda \setminus A}(\varphi_0^2)$$

$$= \zeta^\eta_\Lambda(\mathcal{A} \cap B_\eta = \emptyset) \, \mathrm{E}_{\Lambda \setminus B_\eta}(\varphi_0^2) \,,$$

where the last identity follows from, *e.g.*, the random walk representation of the variance. Now, the estimate (28) implies that

$$\zeta^\eta_\Lambda(\mathcal{A} \cap B_\eta = \emptyset) \geq 1 - c |\log \eta|^{-1} \,,$$

and the conclusion follows from the following standard asymptotics for the Green function of a random walk killed as it exits a box of large radius $R$:

$$G_R(0,0) = (\pi \sqrt{\det \mathcal{Q}})^{-1} \log R + O(1) \,.$$

**The mass.** I only discuss the simpler upper bound on the mass. Denoting by $\mathbb{E}^{(n)}_{i,j}$ the expectation for the random walk starting at $i$ and conditioned on $X_n = j$ (provided the probability of the latter event is positive), we can write

$$\mathrm{E}^\eta(\varphi_0 \varphi_i) \geq \sum_{n \geq 0} \mathbb{E}_0 \left[ \exp\left( -c\eta |\log \eta|^{-1/2} |X_{[0,n]}| \,;\, X_n = i \right) \right]$$

$$= \sum_{n \geq 0} \mathbb{P}_0[X_n = i] \, \mathbb{E}^{(n)}_{0,i} \left[ \exp\left( -c\eta |\log \eta|^{-1/2} |X_{[0,n]}| \right) \right]$$

$$\geq \sum_{n \geq 0} \mathbb{P}_0[X_n = i] \, \exp\left( -c\eta |\log \eta|^{-1/2} \mathbb{E}^{(n)}_{0,i} \left[ |X_{[0,n]}| \right] \right)$$

$$\geq \mathbb{P}_0[X_{\bar n} = i] \, \exp\left( -c\eta |\log \eta|^{-1/2} \mathbb{E}^{(\bar n)}_{0,i} \left[ |X_{[0,\bar n]}| \right] \right) \,,$$

where, in the last expression, the sum has been restricted to the single term $\bar n = \bar n(|i|, \eta) \stackrel{\text{def}}{=} \left[ \eta^{-1/2} |\log \eta|^{3/4} |i| \right]$ if $d = 2$, and $\bar n = \left[ \eta^{-1/2} |i| \right]$ if $d \geq 3$.

The first factor can be bounded by

$$\mathbb{P}_0[X_{\bar n} = i] \geq \frac{c}{\bar n^{d/2}} \, \exp\left[ -c' \frac{|i|^2}{\bar n} \right] \,,$$

for some $c, c' > 0$. It only remains to control the tied-down expectation of the range. For $d \geq 3$, this is trivial, since the bound $\mathbb{E}^{(\bar n)}_{0,i}\left[|X_{[0,\bar n]}|\right] \leq \bar n + 1$ suffices to prove the claim. In dimension 2, in order to get the correct logarithmic correction, more care has to be taken, and one has to resort to the following, less elementary bound, reflecting the effect of the recurrent behavior of the walk,

$$\mathbb{E}^{(A|i|)}_{0,i} \left[ |X_{[0,A|i|]}| \right] \leq c \frac{A|i|}{\log A} \,.$$

I refer to [20, Proposition C.1.] for a proof.

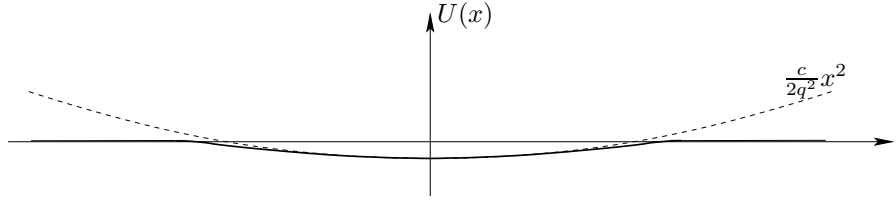

FIG 10. *The pinning potential considered in [46] and its quadratic approximation yielding the effective mass.*

### 5.5. Additional results

**Random potential.** There have also been several works on the study of localization by a random pinning potential, in particular in dimension 1. For example, in [6, 82], it is proved (in a rather general 1-dimensional setup) that if the pinning potential at site $i$ is given by $w + W_i$, with $w \in \mathbb{R}$ a constant, and $W_i$ a family of i.i.d. real-valued random variables with mean zero, then the interface still almost surely gets localized (in the sense that there is a density of pinned sites) even when $w$ is slightly negative (that is, in average the reward is actually a penalty). In [64], it is proved for the same type of models that the presence of disorder induces a smoothening of the phase transition in the sense that it becomes higher order than in the corresponding deterministic case ($W_i \equiv 0$).

In a different spirit, the case of a diluted pinning potential (that is, a pinning potential taking value $\zeta > 0$ or $0$ at each site of the box) is considered in [72]. It is proved that, for $d = 1$ or $2$, the interface is localized (again, in the sense that there is a density of pinned sites) if and only if the sites at which the pinning potential is non-zero have positive density (notice that the disorder is *fixed*, not sampled from some given distribution).

Finally, pinning of an SOS interface by spatial disorder not restricted to a plane but present everywhere in space was studied in [22, 23]. The main results are that: 1) In dimensions $d \geq 3$, the interface is rigid, provided that $\beta$ be large enough and the disorder sufficiently weakly coupled to the field. 2) In dimensions $d \leq 2$, the interface is never rigid.

**Mean-field regime.** The critical behavior of the covariance has also been obtained in a mean-field regime[6] in [46, 47], see also [78]. I briefly describe the setting and the result in order to show the difference with the regime discussed in this section. The measure considered in [46] is the following perturbation of the Gaussian model

$$\mathrm{P}_\Lambda^U(\mathrm{d}\varphi) \propto \prod_{i \in \Lambda} e^{-U(\varphi_i)}\, \mathrm{P}_\Lambda(\mathrm{d}\varphi)\,,$$

---

[6]The term mean-field refers to the fact that in this regime fluctuations are sufficiently small for the interface to remain in the vicinity of the quadratic minimum of the potential ("fluctuations are negligible").

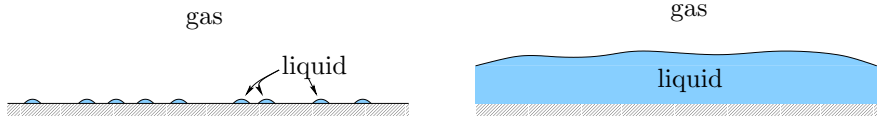

FIG 11. *The wetting transition. Left: partial wetting regime; the liquid phase adsorbed on the wall forms microscopic droplets (of size not diverging with the size of the system). Right: complete wetting regime; the liquid phase adsorbed at the wall forms a mesoscopic layer (its size diverges with the system size). The wetting transition is the phase transition between these two regimes.*

where the self-potential is given by (see Fig. 10)

$$U(x) = -c(e^{-\frac{x^2}{2q^2}} - 1).$$

Then, provided[7] $K\log(1+c^{-1}) < \sqrt{q}$ for some sufficiently large constant $K$ and $0 < c \leq 1$, it is proved that

$$\mathrm{cov}_{\mathbf{P}^U}(\varphi_i, \varphi_j) \leq K\log(q/\sqrt{c})e^{-D\frac{\sqrt{c}}{q}\|j-i\|_2},$$

with the constant $D \to 1$ if $c$ is fixed and $q \to \infty$. The heuristic behind this result is rather clear. Under the above assumption, the quadratic approximation $U(x) = \frac{c}{2q^2}x^2$ holds over a huge range of values of $x$. Over this range of values the measure $\mathbf{P}^U$ behaves like a massive Gaussian model with mass $\mu = \sqrt{c}/q$, and therefore, provided the interface stays mostly there, the exponential decay should be given by this mass. The main part of the proof in [46, 78] was then to control perturbatively this approximation.

The $\delta$-pinning corresponds to an opposite regime, where instead of having a very wide and shallow potential well, one has a very narrow and deep one. It is far less clear *a priori* what the behavior of the correlation lengths should be in this case, since the latter cannot be read from the self-potential.

## 6. Wetting

### 6.1. Description of the model

This section deals with the wetting transition. This is a surface phenomenon of major theoretical and practical interest, and still the object of active study. It occurs each time some substance occupies the bulk of a system, while another substance (or another phase of the same substance) is favoured by the boundary. In this case a layer of the preferred substance can exist in the vicinity of the walls, and the question is to understand its behavior. It turns out that varying some external parameters, say the temperature, can give rise to a phase transition:

---

[7] It is emphasized in [46] that this condition is actually too strong and that the result should be true under the weaker condition that $K\log(1+c^{-1}) < q$, which characterizes the mean-field regime.

In one regime (called *partial wetting*), the thickness of the layer is microscopic[8], while it is mesoscopic in the other regime (called *complete wetting*); see Fig. 11. Practically the manifestation is a transition from a situation where the wall is covered by a multitude of microscopic droplets, to a situation where the wall is covered by a mesoscopic, homogeneous film.

The first rigorous treatment of the wetting transition in the framework of Statistical Mechanics is due to D. Abraham, who, through exact computations, was able to prove its existence for the two-dimensional Ising model [1]. Immediately following this work, there was a series of papers [31, 28, 89, 33] (see also [50]) analyzing the same question in the simpler settings of one-dimensional effective interface models (still through exact computations), and [32] establishing the existence of the wetting transition in the SOS model in any dimension. Only recently has the rigorous analysis of this problem been reconsidered, both providing stronger and more general results in one-dimension, and establishing some preliminary results in higher dimensions. It should be noted that if the problem is well understood, even rigorously, in the two-dimensional Ising model, the situation in the three-dimensional Ising model is still controversial, even among physicists. In particular, it has been suggested [51] that another kind of effective models should be considered if one wants to correctly predict the behavior observed in this model. From a mathematical point of view, these new effective models are too complicated, and even our understanding of the original ones is very far from being satisfactory.

A natural model for the rigorous study of wetting in the framework of effective interface models is obtained by introducing simultaneously a hard wall and a local pinning potential. The latter models the affinity of the wall toward the two equilibrium phases separated by the interface. The basic measure hence takes the form
$$P_\Lambda^{\eta,+}(\,\cdot\,) \stackrel{\text{def}}{=} P_\Lambda^\eta(\,\cdot\,|\,\Omega_+)\,, \tag{32}$$
where $\Omega_+ \stackrel{\text{def}}{=} \{\varphi_i \geq 0, \forall i \in \mathbb{Z}^d\}$.

### 6.2. Main results

Apart from the one-dimensional case, the only rigorous result for discrete models I am aware of is [32], where it is proved for the SOS model in any dimension that there is a phase transition from a regime (at small $\eta$) where there is a density of heights taking value 0, to a regime (at large $\eta$) where this density vanishes.

Actually, there should not be any major obstacle to establishing a complete picture at large values of $\beta$, using a suitable version of the Pirogov-Sinai theory. In fact, there is a claim of such a proof in [26], but the statement seems unlikely (an infinite sequence of first order transitions) and no proof is given there.

---

[8]Microscopic means here that the layer has an average width which is bounded uniformly in the system size, while mesoscopic means that the width diverges (sub-linearly) with the system size.

In the rest of this section, we'll stick to continuous effective interface models. For ease of notations, only the nearest neighbor case will be discussed, but all that is said here can be straightforwardly extended to any finite-range interactions, and probably to long-range interaction with more care.

Since the presence of a hard wall has a strongly repulsive effect on the interface, while the action of a local pinning potential has a strongly localizing effect, their simultaneous presence gives rise to a delicate competition.

The first natural question is whether a non-trivial transition from a localized to a delocalized regime can occur as the intensity of the pinning potential is varied. A convenient way to study this problem at the thermodynamic level is through the *density of pinned sites*, which is defined by

$$\rho_N \stackrel{\text{def}}{=} |\Lambda_N|^{-1} \sum_{i \in \Lambda_N} \mathbf{1}_{\{\varphi_i = 0\}},$$

and its limit $\rho \stackrel{\text{def}}{=} \lim_{N \to \infty} \rho_N$; the existence of this limit follows easily from the FKG property satisfied by the set pinned sites, see [20, Lemma 2.1]. As $\rho_N = |\Lambda_N|^{-1} \eta \frac{d}{d\eta} \log Z^\eta_{\Lambda_N}$, we immediately obtain that $\rho$ is a non-decreasing function of $\eta$. Obviously $\rho(0) = 0$. It is also easy to see that $\rho(\eta) > 0$ for large enough $\eta$. Indeed, since

$$|\Lambda_N|^{-1} \log \frac{Z^{\eta,+}_{\Lambda_N}}{Z^{0,+}_{\Lambda_N}} = \int_0^\eta \frac{1}{\bar\eta} \rho_N(\bar\eta) \, d\bar\eta, \tag{33}$$

the result follows from $Z^{\eta,+}_{\Lambda_N} \geq \eta^{|\Lambda_N|}$ and the existence of a constant $C$ such that $Z^{0,+}_{\Lambda_N} \leq C^{|\Lambda_N|}$. To prove the latter inequality, one can consider a shortest self-avoiding path $\omega$ on $\mathbb{Z}^d$ starting at some site on the boundary of $\Lambda_N$ and containing all the sites of $\Lambda_N$, and use $H_{\Lambda_N}(\varphi) \geq \sum_{n=1}^{|\omega|-1} V(\varphi_{\omega_n} - \varphi_{\omega_{n+1}})/4d$ to write

$$Z^{0,+}_{\Lambda_N} \leq \int \prod_n e^{-V(\varphi_{\omega_n} - \varphi_{\omega_{n+1}})/4d} \leq \left( \int_{-\infty}^\infty e^{-V(x)/4d} dx \right)^{|\Lambda_N|}.$$

This implies that the following critical value is well defined[9]

$$\eta_c \stackrel{\text{def}}{=} \inf\{\eta \,:\, \rho(\eta) > 0\}. \tag{34}$$

The basic question is therefore whether $\eta_c > 0$. The answer turns out to depend both on the dimension and on the tail of the gradient interaction.

1. $V(x) = x^2$: $\eta_c > 0$ if and only if $d \leq 2$.

---

[9]Notice that, since $\rho_N = |\Lambda_N|^{-1} \eta \frac{d}{d\eta} \log Z^\eta_{\Lambda_N}$, the critical value can also be defined as the largest value of $\eta$ for which the free energy of this model is the same as that of the free interface. This is the standard definition in the Physics literature (namely the partial wetting regime is characterized by the fact that the surface tension of the interface along the wall is strictly smaller than that of an interface with the same orientation, but located in the bulk of the system).

2. *V* Lipschitz: $\eta_c > 0$ for all $d \geq 1$.

The "only if" part of the first statement is proved in [19], while the "if" part and the second statement are proved in [29]; heuristic for these results are provided in Subsection 6.3.

**Remark 20.** *The usual derivation of Gaussian effective interface models from a really microscopic model (e.g. a 3-dimensional Ising model above the roughening temperature), relies mainly on a second-order expansion of the free energy along the interface. Such an approximation can be relevant only when large values of gradients have no influence on the physics of the problem. The striking difference in behavior between a harmonic interaction and a Lipschitz one (which may coincide on an arbitrarily large, but finite range of values) shows that it will be very delicate to determine the correct effective interaction (if any) to mimic the behavior of the real system.*

*This of course also shows that results obtained for the harmonic model, depending on the properties of large local fluctuations (e.g. the spikes in the entropic repulsion phenomenon) might not be generic. It would thus be extremely interesting to have some results for models with other tail behaviors.*

By and large, in spite of quite some work, this is unfortunately all that is known when $d \geq 2$. Considering the very detailed pathwise description available both in the case of pure pinning and that of pure entropic repulsion, one might expect to have at least some pathwise informations for the wetting problem. This turns out to be very difficult. Actually, the only results known to date are [90]

- Weak form of delocalization in the whole complete wetting (*i.e.* delocalized interface) regime ($d = 2$):
$$\lim_{N \to \infty} \mathrm{E}^{\eta,+}_{\Lambda_N}(\varphi_i) = \infty \,,$$
  for any $i \in \mathbb{Z}^d$ and $\eta < \eta_c$.
- Strong form of delocalization deeply inside the complete wetting regime ($d = 2$):
$$\mathrm{E}^{\eta,+}_{\Lambda_N}(\varphi_i) \asymp \log N \,,$$
  for any $i \in \Lambda_{\epsilon N}$, $0 < \epsilon < 1$, and $\eta$ sufficiently small.
- Strong form of localization deeply inside the partial wetting regime (*i.e.* localized interface) ($d = 2$):
$$\lim_{N \to \infty} \mathrm{E}^{\eta,+}_{\Lambda_N}(\varphi_i) \leq C \,,$$
  for all $i \in \mathbb{Z}^d$ and $\eta$ sufficiently large. Moreover, the mass is positive:
$$m_{\mathrm{P}^{\eta,+}}(x) > 0 \,,$$
  for all $x \in \mathbb{S}^1$, provided $\eta$ is taken large enough.

As can be seen, the results are still very limited. Particularly annoying is the total absence of pathwise localization results in the full partial wetting regime (and this might be nontrivial when $d \geq 3$, even *deeply* inside the partial wetting regime).

**Open Problem 9.** *Prove the lacking pathwise estimates associated to the wetting transition.*

Let us remark that the approach used in the study of the pinning potential can of course also be applied here, yielding an expression of the form:

$$\mathrm{E}^{\eta,+}_\Lambda(f) = \sum_{A \subset \Lambda} \zeta^{\eta,+}_\Lambda(A)\, \mathrm{E}^+_{\Lambda \setminus A}(f),$$

where $\zeta^{\eta,+}_\Lambda(A) \stackrel{\mathrm{def}}{=} \eta^{|A|} Z^+_{\Lambda \setminus A}/Z^{\eta,+}_\Lambda$. Unfortunately, this is much less useful than before, for two reasons. First, one has very little control over the distribution of pinned sites; essentially the only result that can be easily obtained is the stochastic domination

$$\zeta^{\eta,+}_\Lambda \preccurlyeq \zeta^{\eta}_\Lambda, \qquad (35)$$

stating the rather intuitive fact that there are less pinned sites in the presence of the wall than in its absence. The second reason why this approach seems less promising is that the random walk representation does not apply to the conditional expectation $\mathrm{E}^+_\Lambda(\varphi_i \varphi_j)$ (and this would not even be the most interesting quantity, since the 2-point function does not coincide with the covariance anymore).

Still, there is one non-trivial consequence that can be extracted from (35): The density of pinned sites (and therefore also the height of the interface, thanks to the entropic repulsion results of Section 3, see [90]) diverges continuously as $\eta \downarrow \eta_\mathrm{c}$, provided $\eta_\mathrm{c} = 0$. This is the case, *e.g.*, for the Gaussian model in dimension $d \geq 3$. This shows that the wetting transition is a continuous transition in that case. It is also known to be continuous in dimension 1 [42], but nothing at all in known in the two-dimensional case, although it is clearly expected that the transition is also continuous.

**Open Problem 10.** *Determine the nature of the phase transition in the two-dimensional model.*

**Open Problem 11.** *Study the critical behavior when the transition is second order.*

In dimension $d = 1$, however, the understanding is pretty much complete. After an initial result [71], restricted to a particular choice of underlying random walk, the following result, valid for essentially arbitrary interaction/underlying random walk was proved in [42]. Suppose that the interaction $V : \mathbb{R} \to \mathbb{R} \cup \infty$ is such that $\exp(-V(\,\cdot\,))$ is continuous, $V(0) < \infty$, $\kappa \stackrel{\mathrm{def}}{=} \int \exp(-V(x))\mathrm{d}x < \infty$, $\int x \exp(-V(x))\mathrm{d}x = 0$, and $\int x^2 \exp(-V(x))\mathrm{d}x < \infty$. Then, for the corresponding nearest-neighbor model:

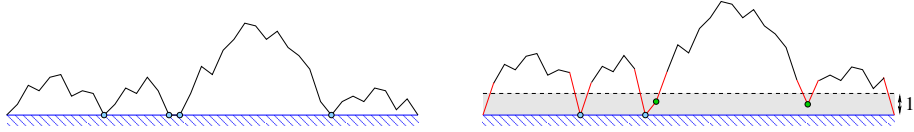

Fig 12. *The transformation: A random subset of the pinned sites is unpinned and the height at the corresponding sites is chosen at random in the interval $(0,1)$. To avoid overcounting, the whole interface, except the pinned sites, is lifted by $1$ before this procedure is applied.*

1. $\eta_c = \kappa/(1 + \sum_{N=1}^{\infty} \kappa^{-N} Z_{\{1,\ldots,N\}}) > 0$.
2. For $\eta \in [0, \eta_c)$, under diffusive scaling, the law of the path converges to that of a Brownian excursion. Moreover, before rescaling, there are only pinned points at a finite (microscopic) distance from the endpoints of the path.
3. For $\eta = \eta_c$, under diffusive scaling, the law of the path converges to that of the absolute value of the Brownian bridge.
4. For $\eta \in (\eta_c, \infty]$, under diffusive scaling, the law of the path converges to the measure concentrated on the constant function taking value zero. Actually, in that case the measure $P^\eta_{\{1,\ldots,N\}}$ converges, without rescaling, to the law of a finitely recurrent, irreducible Markov chain on $\mathbb{R}_+$, which can be explicitly described.

**Remark 21.** *It is worth noting that the corresponding claims also hold when there are free boundary conditions at $N$, see [42].*

### 6.3. Ideas of proofs

**Delocalization at small pinning.** Let us start with the case of $V$ Lipschitz. To simplify, I only consider the special case $V(x) = |x|$.

The idea is to show that the total contribution to the partition function $Z^{\eta,+}_{\Lambda_N}$ of the configurations having a positive density of pinned sites is negligible. This is done by unpinning a random subset of the pinned sites; in this way, one obtains a new set of configurations, whose contribution to $Z^{\eta,+}_{\Lambda_N}$ is exponentially (in $N^d$) larger than the original ones. Let $\epsilon > 0$ and denote by $Z^{\eta,+}_{\Lambda_N}(\epsilon)$ the total weight of the configurations having at least $\epsilon N^d$ pinned sites. Let also $\bar{Z}^{\eta,+}_{\Lambda_N}(\epsilon)$ be the total weight of the configurations having at least $\epsilon N^d$ sites such that $\varphi_i \leq 1$. Starting from any configuration contributing to $Z^{\eta,+}_{\Lambda_N}(\epsilon)$, one can construct a family of configurations contributing to $\bar{Z}^{\eta,+}_{\Lambda_N}(\epsilon)$ by first lifting all the unpinned sites by 1, and then choosing a random subset of the pinned sites and unpinning them, selecting for each of them a height in the interval $(0,1)$ (see Fig. 12). This ensures that the transformation is invertible and thus that there is no overcounting. There are three contributions to the energetical cost of doing that: *(i)* a boundary term, due to the lifting, which is of order $e^{-O(N^{d-1})}$; *(ii)* a factor $1/\eta$ for each unpinned site; *(iii)* a factor never smaller than $e^{-2d}$ for each originally pinned sites, because of the local deformations due to the lifting/unpinning. Therefore, if we denote by $Z_{\Lambda_N}^{\eta,+}(A)$ the total weight of the configurations vanishing exactly on $A$, $\bar{Z}_{\Lambda_N}^{\eta,+}(A,B)$ the total weight of the configurations vanishing exactly on $B$ but such that $\varphi_i \leq 1$ if and only if $i \in A$, and $\bar{Z}_{\Lambda_N}^{\eta,+}(A) = \sum_{B \subset A} \bar{Z}_{\Lambda_N}^{\eta,+}(A,B)$ the total weight of the configurations such that $\varphi_i \leq 1$ if and only if $i \in A$, the above considerations yield

$$\bar{Z}_{\Lambda_N}^{\eta,+}(A,B) \geq e^{-2d|A|} e^{-O(N^{d-1})} \left(\frac{1}{\eta}\right)^{|A|-|B|} Z_{\Lambda_N}^{\eta,+}(A).$$

Therefore

$$\bar{Z}_{\Lambda_N}^{\eta,+}(A) \geq \sum_{B \subset A} \binom{|A|}{|B|} e^{-2d|A|} e^{-O(N^{d-1})} \left(\frac{1}{\eta}\right)^{|A|-|B|} Z_{\Lambda_N}^{\eta,+}(A)$$

$$= e^{-O(N^{d-1})} \left(e^{-2d}(1+\frac{1}{\eta})\right)^{|A|} Z_{\Lambda_N}^{\eta,+}(A).$$

Thus,

$$\frac{Z_{\Lambda_N}^{\eta,+}(\epsilon)}{Z_{\Lambda_N}^{\eta,+}} \leq \sum_{\substack{A \subset \Lambda_n \\ |A| > \epsilon |\Lambda_n|}} \frac{Z_{\Lambda_N}^{\eta,+}(A)}{\bar{Z}_{\Lambda_N}^{\eta,+}(A)} \leq e^{-c\epsilon|\Lambda_n|},$$

for some $c > 0$, since $e^{-2d}(1+\frac{1}{\eta}) > 1$ if $\eta$ is small enough.

The reason this argument does not work immediately in the Gaussian case is that the energetic cost depends in a nonlinear way on the original configuration. However, it is possible (using Jensen inequality) to show that the cost at any site $i$ feeling the lifting is at worst $\exp(\mathrm{E}_{\Lambda_N \setminus A}^+(\varphi_i))$. But, in two dimensions, (21) implies that this cost is finite (since $i$ is then neighboring a pinned site). Therefore, the above argument still applies for the two-dimensional Gaussian model. In higher dimensions, the expectation is infinite so that the argument fails (as it should, since, as we have seen, the interface is localized for any $\eta > 0$ in this case).

**Localization for the Gaussian model with $d \geq 3$.** The proof that the Gaussian interface is always localized when $d \geq 3$ is more delicate. Here, I only give some heuristics.

The main idea is to get, for any $\eta > 0$, a strictly positive lower bound on

$$|\Lambda_N|^{-1} \log \frac{Z_{\Lambda_N}^{\eta,+}}{Z_{\Lambda_N}^{0,+}}.$$

This would then imply that $\rho(\eta) > 0$ for all $\eta > 0$, using (33). It is thus sufficient to construct a suitable set of configurations producing a large enough contribution. The latter configurations are chosen as follows: $\Lambda_N$ is partitionned into cells of spacing $\Delta > 0$ (chosen large enough); then, inside each cell, exactly one site is pinned in the central cube of sidelength $\Delta/5$ (see Fig. 13). The number of choices for the pinned sites is $\exp\left((d\log\Delta + c)|\Lambda_N|/\Delta^d\right)$. On the exponential scale we are interested in, we can replace $Z_{\Lambda_N}^{0,+}$ by $Z_{\Lambda_N}$ (their ratio is $\exp(O(N^{d-1}))$). Therefore, one has to get a lower bound on

$$\frac{Z_{\Lambda_N}^{+,\eta}}{Z_{\Lambda_N}} = \sum_A \eta^{|A|} \frac{Z_{\Lambda_N\setminus A}^+}{Z_{\Lambda_N\setminus A}} \frac{Z_{\Lambda_N\setminus A}}{Z_{\Lambda_N}},$$

where the sum is restricted on sets of pinned sites as described above. The last ratio is easily seen to be larger than $\exp(-c|A|)$, for some $c > 0$, for example using (29) and (30). One is therefore left with estimating the probability that the field pinned on a such a set $A$ (and at the boundary) satisfies the positivity constraint. Due to the fact that the spikes in dimensions 3 and larger are very thin, the cost due to these additional pinned sites is not too large (this is not true when $d \leq 2$, and that is why the interface delocalizes at small pinning: it is entropically too expensive to touch the wall). A careful analysis [19] shows that this probability is bounded below by

$$\frac{Z_{\Lambda_N\setminus A}^+}{Z_{\Lambda_N\setminus A}} \geq \exp\left(-(d\log\Delta - c\log\log\Delta)|\Lambda_N|/\Delta^d\right).$$

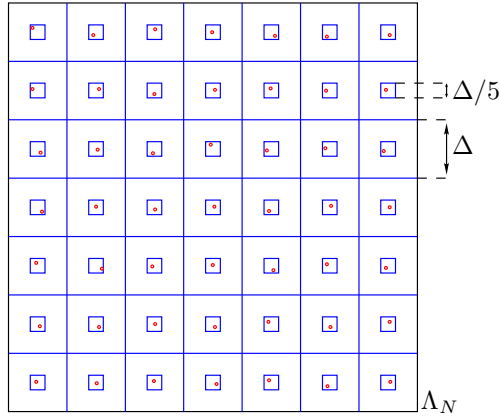

FIG 13. *The box $\Lambda_N$ is decomposed into cells of sidelength $\Delta$, in each of which a smaller cube of sidelength $\Delta/5$ is centered. The partition function is estimated by considering only configurations containing exactly one pinned point inside each of the smaller cubes.*

Notice that this is a result that even goes *beyond* leading order! Putting all this together, we see that $Z_{\Lambda_N}^{+,\eta}/Z_{\Lambda_N}$ (and thus $Z_{\Lambda_N}^{\eta,+}/Z_{\Lambda_N}^{0,+}$) is exponentially large in $|\Lambda_N|$ as soon as $\Delta$ is chosen large enough, depending on $\eta$.

### 6.4. Additional results

**Mean-field regime.** In the already-mentionned papers [46, 47, 78], a model of the wetting transition completely similar to the one discussed for pinning in

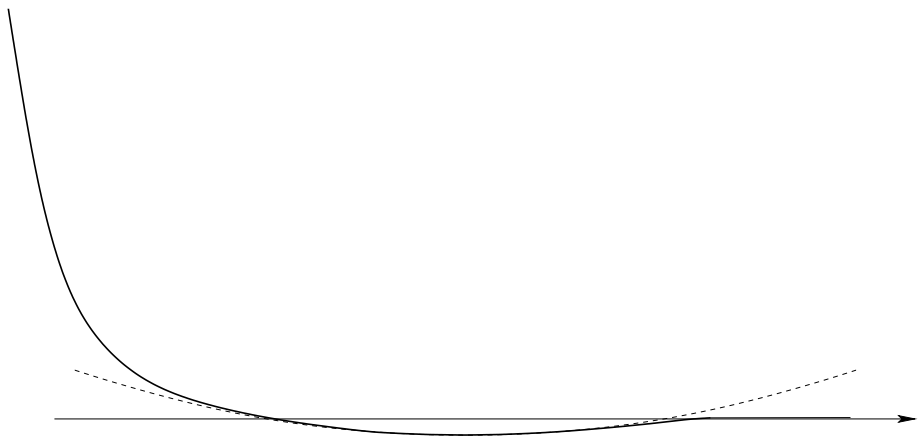

Fig 14. *A sketch of the type of potential considered in [46, 78]. The hard wall is replaced by a fast growth, the pinning part is very shallow, very wide, and very well approximated by a parabola, the curvature of which yields the effective mass of the system.*

Subsection 5.5 has also been studied, and similar results have been obtained. Namely, modelling the attractive potential by a very flat, essentially quadratic potential well (see Fig. 14), it is possible to prove that the mass of the system is again given by the curvature at the bottom of the well.

**Long-range surface/wall interaction.** In the model discussed in this section, the attractive potential has a very short range: It is only in the direct vicinity of the wall that the interface gets a reward. It is claimed in the physics literature that a long-range wall/interface interaction, decaying algebraically with the distance, would yield a first-order wetting transition, or even entirely suppress the transition depending on the decay exponent. I am not aware of any rigorous results in that direction, although physicists might have some exact computations in dimension 1...

**Open Problem 12.** *Analyze the wetting transition in the presence of a long-range wall/interface interaction.*

**Disordered wall.** The model in the presence of a disordered substrate has also been studied in several works, see *e.g.* [35, 36, 49, 21]. Two types of disorder have been considered: random pinning potentials (similar to what is discussed in Subsection 5.5), and rough walls (similar to what is discussed in Subsection 3.3). The main concern of these work is to understand how the surface tension of the system is affected by the disorder, in particular in relation with phenomenological formulas given by physicists (Cassie's law and Wenzel's law).

The works [6, 64], which were already cited in Subsection 5.5, also apply in the case of wetting over a disordered substrate, and the same results as described there also hold in this case.

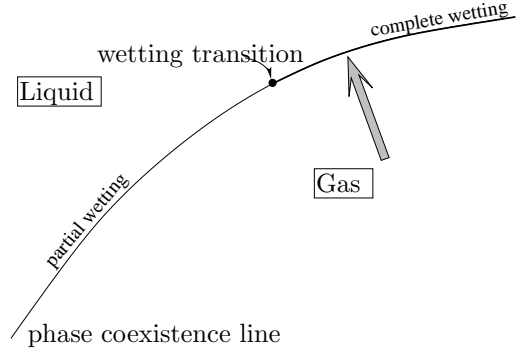

FIG 15. *How does the width of the layer of unstable phase diverges as the coexistence line is approached, in the regime of parameters where complete wetting occurs?*

## 7. Prewetting

### 7.1. Description of the model

The wetting transition occurs exactly on the phase coexistence line: Both the phase occupying the bulk of the system and the phase occupying its boundary (in the complete wetting regime) are thermodynamically stable. One might wonder what happens if one approaches the region of complete wetting from outside the phase coexistence line, see Fig. 15. Indeed, in that case the phase occupying the boundary is not thermodynamically stable anymore, and can only exist sufficiently close to the boundary thanks to the stabilizing effect of the wall. Consequently, even when the film has "infinite" thickness on the phase coexistence line, the latter must be finite away from coexistence. It is then natural to study the way this thickness diverges as the system gets closer and closer to phase coexistence.

In order to model the thermodynamical instability of the wetting layer, we consider the following modification of the measure $P_{\Lambda_N}^{\eta,+}$,

$$P_{\Lambda_N}^{\eta,\lambda,+}(d\varphi) \stackrel{\text{def}}{=} \frac{1}{Z_{\Lambda_N}^{\eta,\lambda,+}} e^{-\lambda \sum_{i \in \Lambda_N} W(\varphi_i)} P_{\Lambda_N}^{\eta,+}(d\varphi), \qquad (36)$$

where $\lambda > 0$ measures the distance to phase coexistence, and the self-potential $W$ should be thought of as being given by $W(x) = |x|$. Actually, since this does not make any essential differences at the level of proofs, at least in the case of continuous heights, one allows any potential $W$ which is convex, increasing on $\mathbb{R}^+$, and satisfies the following growth condition: There exists $f : \mathbb{R}_+ \to \mathbb{R}_+$ such that, for any $\alpha > 0$, we have

$$\limsup_{x \to \infty} \frac{W(\alpha x)}{W(x)} \leq f(\alpha) < \infty.$$

In particular any convex increasing polynomial function is admissible, including our basic example $W(x) = x$. In the latter case, this measure gives a penalization proportional to the volume of the wetting layer (*i.e.* the volume between the wall and the interface modeled by $\varphi$), and $\lambda$ can be interpreted as the difference in free energy between the unstable phase below the interface and the stable phase above.

### 7.2. Main results

It turns out that the actual behavior depends strongly on the nature of the system. When the system is below its roughening transition (e.g. for low temperature Ising model in dimensions $d \geq 3$, or for low temperature discrete effective interface models in dimensions $d \geq 2$), the divergence occurs through an infinite sequence of first-order phase transitions, the *layering transitions*, at which the thickness of the film increases by one microscopic unit. This phenomenon has undergone a detailed rigorous study, at sufficiently large $\beta$, for the discrete SOS model with $W(x) = |x|$ in [44, 30, 77].

Let us now turn to rough interfaces. For systems above their roughening temperatures (e.g. Ising model in dimension $d = 3$ between $T_r$ and $T_c$, the 2-dimensional Ising model at any subcritical temperature, or continuous effective interface models at any temperature), the divergence of the film thickness occurs continuously; this is the so-called *critical prewetting*.

In order to have delocalization of the interface when $\lambda \searrow 0$, it is necessary that $\eta$ be such that the model at $\lambda = 0$ is in the complete wetting regime. Let us therefore introduce the following set of admissible values of $\eta$: $\mathsf{CW} \stackrel{\text{def}}{=} \{0\} \cup \{0 < \eta < \eta_c\}$ (see (34)).

I first discuss the case $d \geq 2$. Set $\mathcal{H}_{\lambda,2} = |\log \lambda|$, and $\mathcal{H}_{\lambda,d} = |\log \lambda|^{1/2}$ when $d \geq 3$. It is shown in [90] that, for all $\eta \in \mathsf{CW}$, uniformly in $\lambda$ small enough,

$$\mathrm{E}^{\eta,\lambda,+}(\varphi_0) \asymp \mathcal{H}_{\lambda,d}. \tag{37}$$

**Open Problem 13.** *Prove that the model is massive, and study the asymptotic behaviour of the mass as $\lambda \to 0$.*

In the one-dimensional case, substantially stronger informations can be extracted, valid for a much larger class of models. This is particularly valuable since we'd like to understand the degree of universality of the critical behavior displayed here. This is the content of the work [68], in which the following model is studied: Configurations of the interface are given by $\varphi \in \mathbb{Z}_+^{\Lambda_N}$, $\Lambda_N = \{-N, \ldots, N\}$, *i.e.* discrete effective interface models are considered; the analysis could also be performed in the case of continuous heights, but would be slightly more cumbersome. Of course, in one dimension, both continuous and

Y. Velenik/Localization and delocalization of random interfaces    161discrete heights models have the same behavior. The probability measure on the set of configurations is given by

$$P_{\Lambda_N}^{0,\lambda,+}(\varphi) \stackrel{\text{def}}{=} \frac{1}{Z_{\Lambda_N}^{0,\lambda,+}} \exp\{-\lambda \sum_{i \in \Lambda_N} W(\varphi_i)\} \prod_{i=-N-1}^{N} \pi(\varphi_{i+1} - \varphi_i),$$

with the usual boundary conditions $\varphi_{-N-1} = \varphi_{N+1} = 0$. As above $\lambda > 0$ and $W$ is a convex function satisfying the same growth condition. $\pi(\,\cdot\,)$ are transition probabilities of an aperiodic one-dimensional integer-valued random walk with increments of zero mean and finite variance. These assumptions on $\pi$ should be optimal (in the sense that otherwise the critical behavior would be different).

**Remark 22.** *The only case previously rigorously studied in the literature is the one-dimensional continuous SOS model with $W(x) = |x|$ [5]. This case turns out to be exactly solvable; in particular, the authors also obtain some explicit constants, which are out of reach in the general case. See also the heuristic discussion in [50].*

**Remark 23.** *Notice that there is no pinning potential (i.e. $\eta = 0$), since this would not change anything in the results (as long as $\eta \in \mathsf{CW}$), but would complicate substantially the analysis. Of course, the analogue of the result stated above for dimensions $d \geq 2$ above can be proved for uniformly strictly convex interactions also in the case $d = 1$ with any $\eta \in \mathsf{CW}$.*

In one dimension, a rather precise description of the full trajectory can be given. Here I don't state the strongest results that can be obtained, but focus on a few specific ones, namely estimates on typical heights, and decay of correlations; see [68] for additional informations.

The typical width of the interface in dimensions $d \geq 2$ does not change qualitatively as the form of the self-potential $W$ is modified. In dimension 1, this is no longer the case, and the quantity describing the typical width of the interface is $\mathcal{H}_\lambda$, which is the unique solution of the equation

$$\lambda \mathcal{H}_\lambda^2 W(2\mathcal{H}_\lambda) = 1\,.$$

In particular, in the case we are mostly interested in, which is $W(x) = x$, $\mathcal{H}_\lambda$ is of order to $\lambda^{-1/3}$.

The first result is analogous to what was obtained in higher dimensions: There exist $c > 0$ and $\lambda_0 > 0$ such that, uniformly in $\lambda < \lambda_0$ and $N > c\mathcal{H}_\lambda^2$,

$$E_{\Lambda_N}^{0,\lambda,+}(\varphi_0) \asymp \mathcal{H}_\lambda\,.$$

Actually, in one dimension, one can obtain more precise estimates on the law of fluctuations on the scale $\mathcal{H}_\lambda$: For any $T$ large enough, and $N$ large compared to $\mathcal{H}_\lambda^2$, there exist $0 < c_2 < c_1 < \infty$ such that

$$\frac{1}{c_1}e^{-c_1 T^{3/2}} \leq P_{\Lambda_N}^{0,\lambda,+}\Big(\varphi_0 \geq T\mathcal{H}(\lambda)\Big) \leq \frac{1}{c_2}e^{-c_2 T^{3/2}}\,,$$

for all $\lambda \in (0, \lambda_0(T)]$ (uniformly on compact subsets).

**Remark 24.** *For a fixed, small $\lambda > 0$, this does not provide informations on arbitrarily large fluctuations on the scale $\mathcal{H}(\lambda)$. It turns out that the behavior deeper in the tail is not universal anymore. In the case of transition probabilities $\pi(\,\cdot\,)$ with Gaussian-like tails, the same result holds uniformly in all small $\lambda$, for all $T$ large enough. However, when the tail of $\pi(\,\cdot\,)$ becomes fatter, the behavior changes qualitatively. It is therefore rather remarkable that one can still extract some universal information (the exponent $3/2$).*

Moreover, in one dimension, it is possible to obtain estimates on the decay of correlations: There exists $c, c_1, c_2 > 0$, $\lambda_0 > 0$ and $\delta > 0$ such that, uniformly in $\lambda < \lambda_0$ and $N > c\mathcal{H}_\lambda^2$,

$$\frac{\mathcal{H}_\lambda^2}{c_1} \exp\left(-c_1 \mathcal{H}_\lambda^{-2} \|j-i\|_2\right) \leq \mathrm{cov}_{\mathrm{P}_{\Lambda_N}^{0,\lambda,+}}(\varphi_i, \varphi_j) \leq \frac{\mathcal{H}_\lambda^{5/2}}{c_2} \exp\left(-c_2 \mathcal{H}_\lambda^{-2} \|j-i\|_2\right).$$

In particular, the mass satisfies $m_{\mathrm{P}^{0,\lambda,+}} \asymp \mathcal{H}_\lambda^{-2}$. (The lower bound on the covariance, and therefore the upper bound on the mass are not stated in [68], but are very easily obtained; a sketch of the argument is given below.)

### 7.3. Some heuristic

Even though the proofs are somehow technical, it is not difficult to get some intuition about the results, in particular the critical exponents. I only discuss the case $W(x) = |x|$, but the general case is similar.

Suppose $d = 1$. First, if one expects that the interface does indeed remain at an average distance $H$ from the wall, then the corresponding energetic cost is of order $\lambda H N$. Once this energetic contribution is removed, one is left with a pure random walk problem: what is the entropic cost for a random walk conditioned to stay positive, to remain below some fixed level $H$ for a time $N$? This is the same problem which is discussed in Section 4, where we have seen that the probability of such an event is of order $e^{-O(N/H^2)}$. Therefore an energy/entropy balance gives, $\lambda H N \sim N/H^2$, i.e. $H \sim \lambda^{-1/3}$.

The same heuristic also applies when $d \geq 2$, provided one uses the corresponding results of Section 4, and yields the claim (37).

Let us return to $d = 1$. To estimate the covariance, we use the standard duplication trick. Namely let $\varphi$ and $\varphi'$ be two independent copies of the process. Then $\mathrm{cov}_{\mathrm{P}_{\Lambda_N}^{0,\lambda,+}}(\varphi_i, \varphi_j) = \frac{1}{2}\mathrm{E}_{\mathrm{P}_{\Lambda_N}^{0,\lambda,+} \otimes \mathrm{P}_{\Lambda_N}^{0,\lambda,+}}[(\varphi_i - \varphi'_i)(\varphi_j - \varphi'_j)]$. Now, suppose that the two paths really never leave the tube between the wall and height $\mathcal{H}_\lambda$, and assume one can neglect the energetic term, then these two paths will meet in a time of order $\mathcal{H}_\lambda^2$, and therefore an elementary coupling between the paths would yield the upper bound on the covariance, since if they meet the expectation in the RHS is necessarily zero by symmetry (under an interchange of the two paths after their first common point, the product changes sign).

To prove the lower bound, it suffices to force $\varphi_k$ to stay between $0$ and $\mathcal{H}_\lambda$, while $\varphi'_k$ stays between $2\mathcal{H}_\lambda$ and $3\mathcal{H}_\lambda$, for all $k \in \{i, \ldots, j\}$. The cost of doing that is of order $e^{-O(\|j-i\|_2/\mathcal{H}_\lambda^2)}$, and under this event, the expectation in the RHS is bounded below by $O(H_\lambda^2)$.

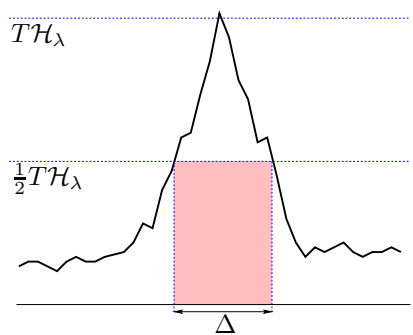

Fig 16. *Idea behind the proof of the tail exponent.*

Finally, the "tail" exponent $3/2$ can also be understood easily (at this heuristic level); see Fig. 16. If one assumes that the width of the excursion to height $T\mathcal{H}_\lambda$, *i.e.* the size of the interval where the path is actually above, say, $\frac{1}{2}T\mathcal{H}_\lambda$ is $\Delta$, then we have an energetically cost $\lambda\frac{1}{2}T\mathcal{H}_\lambda\Delta$ and an entropic cost (assuming Gaussian tails) $e^{-O(T^2\mathcal{H}_\lambda^2/\Delta)}$, which once balanced give $\Delta \sim \sqrt{T\mathcal{H}_\lambda/\lambda}$ and therefore the probability of such a deviation should be, using $\mathcal{H}_\lambda \sim \lambda^{-1/3}$, of order $e^{-O(T^{3/2})}$.

## Acknowledgments

These notes have originally been written for a 6-hours mini-course given in the workshop *Topics in Random Interfaces and Directed Polymers*, Leipzig, September 12-17, 2005. I would like to express my gratitude to the organizers for inviting me to give this mini-course. I also thank E. Bolthausen, J.-D. Deuschel, F. Dunlop, G. Giacomin and H. Sakagawa for discussions and for pointing out errors and omissions in the first draft of these notes. The careful reading and the comments of the anonymous referee have also been very useful.